\documentclass[a4paper,11pt,DIV=13,headinclude]{scrartcl}

 \usepackage{scrlayer-scrpage}
  \clearpairofpagestyles
  \ihead{Disarlo, Pan, Randecker, Tang}
  \ohead{Large-scale geometry of the saddle connection graph}
 \KOMAoptions{headsepline = true}

\usepackage{cmap}           % to search and copy ligatures
\usepackage[utf8]{inputenc} % for Linux computer and Mac
\usepackage[T1]{fontenc}    % to search for ligatures in the pdf

  \usepackage{amsmath,amssymb,amsthm,amsfonts,dsfont,wasysym,pifont,xspace,xcolor}
  \usepackage{url}
  \usepackage{enumitem}
  \usepackage{mathrsfs}
  \usepackage{mathtools}	% für \coloneqq
  \usepackage{xfrac}        % for \sfrac{}{}
  \usepackage[
     breaklinks,colorlinks,
     citecolor=blue,linkcolor=blue,urlcolor=teal
  ]{hyperref}
  \usepackage[nameinlink,capitalise,noabbrev]{cleveref}
  \usepackage{tikz}
  \usetikzlibrary{
    arrows,
    calc,
    matrix,
    decorations.pathmorphing,
    decorations.markings
    }

  \begingroup
    \makeatletter
    \@for\theoremstyle:=definition,remark,plain\do{%
      \expandafter\g@addto@macro\csname th@\theoremstyle\endcsname{%
        \addtolength\thm@preskip\parskip
        }%
      }
  \endgroup

  \makeatletter
    \def\tagform@#1{\maketag@@@{%
     \textbf{(\ignorespaces#1\unskip\@@italiccorr)}}}%
     \renewcommand{\eqref}[1]{\textup{\maketag@@@{(\ignorespaces%
          {\ref{#1}}\unskip\@@italiccorr)}}}
  \makeatother
  
  \theoremstyle{plain}
  \newtheorem{theorem}[equation]{Theorem}
  \newtheorem{proposition}[equation]{Proposition}
  \newtheorem{corollary}[equation]{Corollary}
  \newtheorem{lemma}[equation]{Lemma}

  \theoremstyle{definition}
  \newtheorem{definition}[equation]{Definition}
  
  \newtheorem{example}[equation]{Example}

  \numberwithin{equation}{section}

  \newcommand{\set}[1]{\ensuremath{\left\{ {#1} \right\}}\xspace}

  \newcommand{\st}{\ensuremath{\,\, \colon \,\,}\xspace}
  \newcommand{\from}{\ensuremath{\colon \thinspace}\xspace}

  \DeclareMathOperator{\diam}{diam}
  \DeclareMathOperator{\Aff}{Aff}

  \DeclareMathOperator{\Aut}{Aut}

  \DeclareMathOperator{\QI}{QI}

  \newcommand{\R}{\ensuremath{\mathbb{R}}\xspace}
  
  \newcommand{\SL}{\ensuremath{\mathrm{SL}(2,\R)}\xspace}
  
  \newcommand{\RP}{\R \mathrm{P}^1}

  \newcommand{\MF}{\ensuremath{\mathcal{MF}}\xspace}
  \newcommand{\PMF}{\ensuremath{\mathcal{PMF}}\xspace}
  \newcommand{\AF}{\ensuremath{\mathcal{AF}}\xspace}
  \newcommand{\C}{\ensuremath{\mathcal{C}}\xspace}
  \newcommand{\A}{\ensuremath{\mathcal{A}}\xspace}
  \renewcommand{\AC}{\ensuremath{\mathcal{AC}}\xspace}
  \newcommand{\G}{\ensuremath{\mathcal{G}}\xspace}
    \newcommand{\I}{\ensuremath{\mathcal{I}}\xspace}
   \newcommand{\straight}{\ensuremath{\mathbf{s}}\xspace}
  \newcommand{\ZZ}{\mathcal{Z}}

\title    {Large-scale geometry \\ of the saddle connection graph}
\author{
Valentina Disarlo
\thanks{Mathematisches Institut, Universit\"at Heidelberg, Im Neuenheimer Feld 205, 69120 Heidelberg, Germany,
{\href{mailto:vdisarlo@mathi.uni-heidelberg.de}{vdisarlo@mathi.uni-heidelberg.de}}}
\and Huiping Pan
\thanks{Department of Mathematics, Jinan University, 601 Huangpu Road, Tianhe, Guangzhou, Guangdong, 510632, China 
{\href{mailto:panhp@jnu.edu.cn}{panhp@jnu.edu.cn}}}
\and Anja Randecker
\thanks{Mathematisches Institut, Universität Heidelberg, Im Neuenheimer Feld 205, 69120 Heidelberg, Germany,
{\href{mailto:randecker@mathi.uni-heidelberg.de}{randecker@mathi.uni-heidelberg.de}}}
\and Robert Tang
\thanks{Department of Pure Mathematics, Xi'an Jiaotong--Liverpool University, 111 Ren'ai Road, Suzhou Industrial Park, Suzhou, Jiangsu, 215123, China,
{\href{mailto:robert.tang@xjtlu.edu.cn}{robert.tang@xjtlu.edu.cn}}}
}
\date     {December 9, 2020}

\begin{document}

\maketitle

\begin{abstract}
 We prove that the saddle connection graph associated to any half-translation surface is $4$--hyperbolic and uniformly quasi-isometric to the regular countably infinite-valent tree. Consequently, the saddle connection graph is not quasi-isometrically rigid. We also characterise its Gromov boundary as the set of straight foliations with no saddle connections.
 In our arguments, we give a generalisation of the unicorn paths in the arc graph which may be of independent interest.
\end{abstract}

\section{Introduction}

The large-scale geometry of combinatorial complexes associated to topological surfaces such as the curve graph and arc graph has played a substantial role in answering important questions on mapping class groups, Teichm\"uller theory, and hyperbolic $3$--manifolds \cite{BBF, masur_minsky, brock_canary_minsky, mosher}. More recently, analogous complexes associated to Euclidean cone metrics have attracted considerable attention; see \cite{minsky_taylor_17, nguyen_17, nguyen_18, disarlo_randecker_tang_18, pan_18, fu_leininger_19}. In this paper, we consider the large-scale geometric features of the \emph{saddle connection graph} $\A(S,q)$ associated to a half-translation surface $(S,q)$; this has the saddle connections on $(S,q)$ as vertices, with edges representing pairs of non-crossing saddle connections.

\begin{theorem}[Saddle connection graph is a quasitree]\label{thm:main} For any half-translation surface $(S,q)$, there exists a surjective $(30,29)$--quasi-isometry from the saddle connection graph $\A(S,q)$ to the regular (countably) infinite-valent tree $T_\infty$.
\end{theorem}

Consequently, \emph{all} saddle connection graphs are \emph{uniformly} quasi-isometric to one another.
This demonstrates a disparity between the large-scale and fine-scale geometry of $\A(S,q)$. Indeed, in previous work of the authors \cite{disarlo_randecker_tang_18, pan_18}, it is shown that any isomorphism between saddle connection graphs is induced by an affine diffeomorphism between the underlying half-translation surfaces.
Thus, while the isomorphism class of $\A(S,q)$ completely determines the affine equivalence class of $(S,q)$, the quasi-isometry class reveals no information about the underlying geometric structure.

This disparity is further exemplified when we consider the notion of quasi-isometric rigidity.
The \emph{quasi-isometry group}~$\QI(X)$ of a metric space $X$ is its group of quasi-isometries considered up to the relation of finite distance.
We say $X$ is \emph{quasi-isometrically rigid} if the natural map ${\mathrm{Isom}(X) \rightarrow \QI(X)}$ is an isomorphism. On the one hand, $\Aut(\A(S,q))$ is isomorphic to the affine diffeomorphism group $\Aff^\pm(S,q)$ \cite{disarlo_randecker_tang_18, pan_18} which is countable. On the other hand, $\QI(\A(S,q)) \cong \QI(T_\infty)$ is uncountable as it contains an embedded copy of $\Aut(T_\infty)$. Therefore, $\Aut(\A(S,q)) = \mathrm{Isom}(\A(S,q)) \rightarrow \QI(\A(S,q))$ cannot be an isomorphism.

\begin{corollary}[No QI-rigidity]
 The saddle connection graph $\A(S,q)$ is not quasi-isometrically rigid. \qed
\end{corollary}

The above results are somewhat surprising in light of quasi-isometric rigidity of the curve graph \cite{rafi_schleimer_11}, the mapping class group \cite{hamenstaedt_05, behrstock_kleiner_minsky_mosher_12, bowditch_18}, and of Teichm\"uller space \cite{eskin_masur_rafi_18,bowditch_16,bowditch_20}.
 In particular, Rafi and Schleimer showed that two curve graphs are quasi-isometric to one another if and only if they are isomorphic.
It would be interesting to see if quasi-isometric rigidity also holds for the arc graph.

Let us now turn our attention to the boundary structure of $\A(S,q)$. An immediate corollary of \cref{thm:main} is the following.

\begin{corollary}[Totally disconnected Gromov boundary]
The Gromov boundary of the saddle connection graph $\partial \A(S,q)$ is totally disconnected. \qed
\end{corollary}

Since $\A(S,q)$ isometrically embeds into the arc-and-curve graph $\AC(S)$ \cite{minsky_taylor_17}, we obtain a natural embedding $\partial\A(S,q) \hookrightarrow \partial \AC(S)$ of their respective Gromov boundaries. It is well-known that the inclusion of the curve graph $\C(S) \hookrightarrow \AC(S)$ is a quasi-isometry, (see for instance~\cite{schleimer, kork_pap}) and so $\partial \AC(S)$ is homeomorphic to $\partial \C(S)$.
Klarreich naturally identifies the Gromov boundary of the curve graph $\partial \C(S)$ with the space of arational topological foliations on $S$ \cite{klarreich_18}. Therefore, a natural question is to determine which foliations arise in $\partial\A(S,q) \subset \partial \AC(S)$ under this identification.
We show that these are exactly the foliations that are compatible with the Euclidean structure.
(One could also consider the embedding of the saddle connection graph into the arc graph and results from \cite{pho-on_17}, however, we choose to work with $\AC(S)$ as this simplifies some arguments.)

\begin{theorem}[Boundary of the saddle connection graph]
 The boundary of the saddle connection graph $\partial \A(S,q)$ is the set of all straight foliations on $(S,q)$ that do not contain saddle connections.
\end{theorem}

While the large-scale geometry of the saddle connection graph itself reveals no information, the large-scale geometry of the \emph{embedding} $\A(S,q) \hookrightarrow \AC(S)$ does in fact recover the underlying half-translation structure up to $\SL$--deformations. This follows from an argument of \cite{duchin_leininger_rafi_10}, which shows that the $\SL$--orbit of $(S,q)$ can be determined by any pair of uniquely ergodic foliations from $\partial \A(S,q)$.

In a coarse sense, the saddle connection graph can be regarded as much ``smaller'' than the arc-and-curve graph. To be more precise, we can consider the notion of \emph{asymptotic dimension}, a quasi-isometry invariant. By \cite{bell_fujiwara, fujiwara_whyte} a Gromov hyperbolic space $X$ has asymptotic dimension at most $1$ if and only if $X$ is quasi-isometric to a tree. In such cases, $\partial X$ is necessarily totally disconnected. However, since $\partial\C(S)$ is connected (except for sporadic surfaces in low genus)~\cite{gabai, gabai_09}, so is $\partial \AC(S)$ and thus $\AC(S)$ must have asymptotic dimension at least~$2$. (The arc graph also has asymptotic dimension at least $2$ since its boundary contains the space of arational topological foliation \cite{pho-on_17}, which is itself connected.)

In conclusion, our results demonstrate that combinatorial complexes associated to \emph{geometric structures} on surfaces can have differing behaviour to their topological counterparts. Indeed, in the setting of the saddle connection graph, we show that it is the \emph{fine-scale} geometry to which one should play closer attention, rather than the large-scale geometry.

\paragraph{Acknowledgements}
The first author acknowledges support from the Olympia Morata Habilitation Programme of Universit\"at Heidelberg and from the European Research Council under ERC-Consolidator grant
614733 (GEOMETRICSTRUCTURES). The second author acknowledges the support of National Natural Science Foundation of China  NSFC 11901241. The third author acknowledges the support of NYU-ECNU Institute of Mathematical Sciences at NYU Shanghai. The fourth author acknowledges support from a JSPS KAKENHI Grant-in-Aid for Early-Career Scientists (No.~19K14541). Substantial progress on this project was carried out while some of the authors met at NYU Shanghai which we thank for its kind hospitality.
The authors also thank Kasra Rafi and Saul Schleimer for helpful background discussions on arc graphs and curve graphs and Richard Webb for helpful comments.

\section{Background}

In this section, we recall some basic results about Gromov boundaries, half-translation surfaces, foliations, and combinatorial complexes associated to surfaces.

\subsection{Gromov boundaries}

We begin by recalling the definition of a quasi-isometric embedding. For more background on the material presented in this subsection, see Bridson--Haefliger \cite{Br}. 
\begin{definition}[Quasi-isometric embedding]
Let $(X,d_X)$ and $(Y,d_Y)$ be metric spaces. Given $K\geq 1$ and $C \geq 0$, we say that a map $f: X \to Y$ is a $(K, C)$--\emph{quasi-isometric embedding} if
       $$ K^{-1}\cdot d_X(x,x')-C\leq d_Y(f(x),f(x'))\leq K \cdot d_X(x,x')+C$$
for all $x,x'\in X$. If there also exists a constant $M \geq 0$ such that for any $y\in Y$, there exists~$x\in X$ with $d_Y(y,f(x))\leq M$, then $f$ is called a $(K,C)$--\emph{quasi-isometry}.
\end{definition}
Next, we recall the definition of Gromov hyperbolicity via $k$--centres of geodesic triangles.
\begin{definition}[$k$--centre, Gromov hyperbolic space]
Let $(X,d)$ be a metric space and $k \geq 0$.
A \emph{$k$--centre} of a geodesic triangle $T$ in $X$ is a point $x\in X$ that lies within distance $k \geq 0$ of each of the three sides of $T$.
We say $(X,d)$ is \emph{$k$--hyperbolic} if every geodesic triangle has a $k$--centre. Furthermore, we say $(X,d)$ is \emph{(Gromov) hyperbolic} if it is $k$--hyperbolic for some $k \geq 0$.
\end{definition}

Let us now recall the construction of the Gromov boundary of a hyperbolic space.
\begin{definition}[Gromov product]
Let $(X, d)$ be a metric space. For any $x, y, z \in X$, define the \emph{Gromov product} of $x$ and $y$ with respect to $z$ as
$$(x|y)_z := \frac{1}{2}(d(x, z) + d(y, z) - d(x,y)).$$
\end{definition}
In a tree, the Gromov product $(x|y)_z$ measures how long the geodesics from $z$ to $x$ and from~$z$ to~$y$ coincide before diverging. In Gromov hyperbolic spaces, the Gromov product coarsely measures how long two geodesics fellow-travel.

\begin{definition}[Admissible sequence]
Let $(X,d)$ be a Gromov hyperbolic space with a basepoint~$x \in X$. We say that a sequence $x_n$ of points in $X$ is \emph{admissible} if
$$\liminf_{i, j \to \infty} (x_i | x_j )_x = \infty .$$
We call two admissible sequences $(x_n)$ and $(y_n)$ \emph{equivalent}, and write $(x_n) \sim (y_n)$, if
$$ \liminf_{i,j\to \infty} (x_i | y_j )_x  = \infty.$$
\end{definition}
It is not difficult to see that this definition does not depend on the choice of a basepoint, and that $\sim$ does indeed define an equivalence relation.

\begin{definition}[Gromov boundary]
Let $(X, d)$ be a Gromov hyperbolic space with a basepoint~$x\in X$.
The \emph{Gromov boundary} $\partial X$ of $X$ is the set of equivalence classes of admissible sequences, with the topology defined by equipping each point $p \in \partial X$ with a neighbourhood basis comprising the sets
$$U(p, R) :=\{q \in \partial X~|~ \exists (x_n) \in p, (y_n)\in q
\mbox{ such that } \liminf_{i,j \to \infty} (x_i | y_j )_x >  R \}$$
for all $R > 0$. The topology on $\partial X$ does not depend on the choice of basepoint.
\end{definition}

The following is a well-known fact. 

\begin{proposition}[Quasi-isometries induce homeomorphisms on the boundary \cite{Br}]
Let $X$ and~$Y$ be Gromov hyperbolic spaces and $f : X \to Y$ be a quasi-isometric embedding.
Then $f$ induces a topological embedding $\partial f: \partial X \hookrightarrow \partial Y$. If $f$ is a quasi-isometry then $\partial f$ is a homeomorphism. \qed
\end{proposition}

\subsection{Half-translation surfaces}\label{ssec:MinskyTaylor}

Throughout this paper, we consider a closed, connected, compact topological surface $S$ equipped with a finite, non-empty set of marked points $\ZZ \subset S$.

\begin{definition}[Essential arcs and curves]
A  (simple) \emph{arc} on $(S, \ZZ)$ is a map \hbox{$\alpha: [0,1] \to S$} such that $\alpha|_{(0,1)}: (0,1) \to S \setminus \ZZ$ is an embedding and $\alpha(0), \alpha(1) \in \ZZ$. An arc is \emph{essential} if $\alpha$ is not homotopic (with respect to $\ZZ$) to a point.
Similarly, a (simple closed) \emph{curve} on $(S, \ZZ)$ is an embedding $\gamma: S^1 \to S \setminus \ZZ$. A curve is \emph{essential} if no component of~$S \setminus \gamma$ is a disk or a disk containing exactly one point of $\ZZ$.
\end{definition}

Two arcs or curves are said to be \emph{disjoint} if their images have disjoint interiors.
We shall identify an arc or curve with its image and we consider only simple essential arcs and curves unless otherwise specified.

The main objects of interest in this paper are half-translation surfaces and the combinatorics of their saddle connections. We recall here a few basic definitions, see \cite{zorich_06, strebel} for more on half-translation surfaces.

\begin{definition}[Half-translation surface]
A \emph{half-translation surface} $(S,q)$ is defined by a compact connected surface~$S$, a finite set $\ZZ \subset S$, and a half-translation structure on~$S \setminus \ZZ$, i.e.\ a maximal atlas on~$S \setminus \ZZ$ such that the transition maps are of the form $z \mapsto \pm z + c$ with $c \in \mathbb C$. The set $\ZZ$ is called \emph{singularities} of $(S,q)$ and a singularity is called \emph{removable} if the half-translation structure can be extended to it.
\end{definition}

Every half-translation structure on a surface $S$ is equivalent to a complex structure together with a meromorphic quadratic differential~$q$, where $\ZZ$ contains the set of zeros and simple poles~of~$q$.

We equip $(S,q)$ with a locally Euclidean metric by pulling back the Euclidean metric from $\mathbb{C}$.

\begin{definition}[Saddle connection]
A \emph{saddle connection} on $(S,q)$ is a geodesic segment whose endpoints are (not necessarily distinct) singularities and which has no singularities in its interior.
\end{definition}

Whenever we have a half-translation structure on $(S,\ZZ)$, we implicitly assume that the set of singularities of $(S,q)$ is equal to the set $\ZZ$ of marked points. In other words, we are working in the \emph{fully-punctured} setting in the sense of Minsky--Taylor \cite{minsky_taylor_17}.
We shall recall some useful constructions and results from their paper.

Let $X = (S,q) \setminus \ZZ$ be the half-translation surface $(S,q)$ with its singularities deleted, and~$\hat X$ be the metric completion of its universal cover $\tilde X$. The covering map $\tilde X \rightarrow X$ extends to a branched covering $\hat X \rightarrow (S,q)$, where the set of completion points $\hat\ZZ := \hat X \backslash \tilde X$ coincides with the preimage of $\ZZ$. In particular, saddle connections on $(S,q)$ lift to straight segments connecting distinct completion points. Furthermore, $\hat X$ is a CAT(0)--space, and so every path has a unique geodesic representative in its proper homotopy class.
By a \emph{complete geodesic ray} on $\hat X$, we mean either an infinite geodesic ray or a geodesic path that terminates at some point in $\hat\ZZ$.
A \emph{complete geodesic} on $\hat X$ is a geodesic that decomposes into two complete geodesic rays whenever it is cut along any point in its interior.

\begin{theorem}[Minsky--Taylor compactification \cite{minsky_taylor_17}] \label{thm:compactification}
 There exists a $\pi_1(X)$--equivariant compactification $\mathbf{X}$ of $\hat X$ that is homeomorphic to the closed disc and satisfies the following:
\begin{itemize}
 \item $\tilde X$ coincides with the interior of $\mathbf{X}$ and $\hat\ZZ$ is contained in $\partial \mathbf{X} \cong S^1$,
 \item Each complete geodesic ray on $\hat X$ determines a unique endpoint on $\partial \mathbf{X}$ which, in addition, belongs to $\hat \ZZ$ if and only if the given ray has finite length,
 \item Each complete geodesic on $\hat X$ extends to an embedded path in $\mathbf{X}$ with distinct endpoints on~$\partial \mathbf{X}$,
 \item For any pair of distinct points on $\partial \mathbf{X}$, there exists a complete geodesic on $\hat X$ with the given pair of points as endpoints,
 \item Each parabolic element of $\pi_1(X)$ fixes exactly one point on $\partial \mathbf{X}$ which, moreover, belongs to $\hat \ZZ$; whereas each non-parabolic element fixes a pair of distinct points on $\partial \mathbf{X}$, neither of which belong to $\hat \ZZ$.\qed
\end{itemize}
\end{theorem}

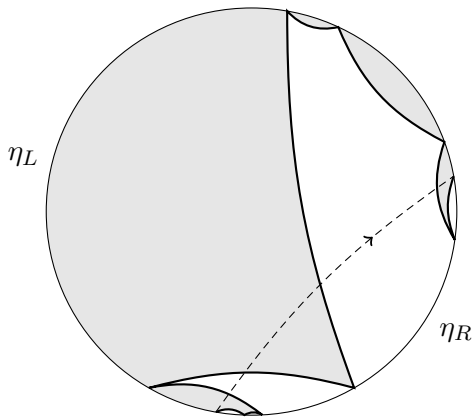
\begin{figure}
 \centering
  \begin{tikzpicture}[scale=0.9, decoration={markings, mark=at position .7  with {\arrow[line width=0.3mm]{>}}}]
   \fill[color=gray!20] (-120:3)
   arc (-120:-100:3) to[bend left=30]
   (-92:3) to[bend left=30]
   (-87:3) to[bend right=25]
   (-120:3);
   \fill[color=gray!20] (-120:3) to[bend left=15]
   (-60:3) to[bend left=10]
   (80:3) arc (80:240:3);
   \fill[color=gray!20] (80:3) to[bend right=30]
   (65:3) to[bend right=20]
   (20:3) to[bend right=25]
   (-8:3) to[bend left=20]
   (10:3) arc (10:80:3);

   \draw (0,0) circle (3cm);

   \draw[densely dashed, postaction={decorate}] (-100:3) to[bend left=10] (10:3);

   \draw[thick] (-100:3) to[bend left=30]
   (-92:3) to[bend left=30]
   (-87:3) to[bend right=25]
   (-120:3) to[bend left=15]
   (-60:3) to[bend left=10]
   (80:3) to[bend right=30]
   (65:3) to[bend right=20]
   (20:3) to[bend right=25]
   (-8:3) to[bend left=20]
   (10:3);

   \draw (170:3) node[above left]{$\eta_L$};
   \draw (-30:3) node[below right]{$\eta_R$};
  \end{tikzpicture}
  \caption{An oriented arc (dashed) and its geodesic representative (bold) in the Minsky--Taylor compactification. The left side of the geodesic is shaded.}
  \label{fig:MT_compactification}
\end{figure}

Using this compactification, Minsky and Taylor define the two \emph{sides} of a complete (oriented) geodesic $l$ on $\hat X$ as follows. The endpoints $\partial l$ separate $\partial \mathbf{X}$ into two open intervals, which we denote~$\eta_L$ and $\eta_R$, that lie to the left and right of $l$ respectively. By the Jordan curve theorem,~$\mathbf{X} \backslash l$ is a disjoint union of at least two connected components. Note that $l$ could possibly pass through points of $\hat\ZZ$ (see \cref{fig:MT_compactification}), and so $\mathbf{X} \backslash l$ could have more than two components.
Consider those components that have non-trivial intersection with~$\eta_L$; we call the closure of their union the \emph{left side} of $l$ and denote it $H_L(l)$. The \emph{right side} $H_R(l)$ is defined similarly.

\begin{lemma}[Convex sides]\label{lem:convex}
 Let $l$ be a complete geodesic on $\hat X$. Then $H_L(l) \cap \hat X$ and $H_R(l) \cap \hat X$ are both convex.
\end{lemma}

\proof
The proof is essentially the same as that of Corollary 2.4 in \cite{minsky_taylor_17}. Let $x$ and $y$ be points in $H_L(l) \cap \hat X$, and let $g$ be the unique geodesic in $\hat X$ connecting them. Suppose for a contradiction that $g$ is not completely contained in $H_L(l)$. Let $g'$ be a maximal subpath of $g$  that is disjoint from $H_L(l)$. Then the endpoints of $g'$, which we shall denote $x'$ and $y'$, must lie on $l \subseteq \partial H_L(l)$. But the subpath of $l$ from $x'$ to $y'$ is also a geodesic, contradicting the unique geodesic property for CAT(0)--spaces.
\endproof

We will later consider the geodesic representatives of arcs and curves.
For our setting, we need to extend the notion of representatives.

Let $\alpha$ first be the homotopy class of an arc on $(S,\ZZ)$, and let $a$ be a representative. We call a path~$a'$ on $S$ a representative of $\alpha$ if it is related to $a$ by a proper homotopy that can be lifted to $\hat X$. In other words, we can homotope $a$ so that it touches marked points, but we cannot homotope it past them. Since $\hat X$ is CAT(0), the proper homotopy class $\alpha$ has a unique geodesic representative~$\alpha_q$ on $(S,q)$; moreover, $\alpha_q$ is a finite concatenation of saddle connections. Note that a path on $(S,q)$ could be the geodesic representative of more than one homotopy class of~arcs.

If $\alpha$ is a curve, there are two possibilities for its geodesic representatives. Either there is a unique geodesic representative $\alpha_q$ which, moreover, is a concatenation of saddle connections; or there is a unique maximal open Euclidean cylinder foliated by closed geodesic loops representing~$\alpha$. In the latter case, write $\alpha_q$ for the set of saddle connections forming the boundary of the cylinder.

The following is a consequence of the convexity of sides.

\begin{proposition}[Geodesic representatives of arcs or curves \cite{minsky_taylor_17}]\label{prop:MT}
Let $\alpha$ and $\beta$ be a pair of disjoint arcs or curves on $(S,\ZZ)$.
Then no saddle connection appearing on $\alpha_q$ transversely intersects a saddle connection appearing on $\alpha_q$ or $\beta_q$. \qed
\end{proposition}

\subsection{Foliations}\label{sec:background_foliations}
In this subsection, we recall some basic properties of measured foliations (for more background, see Fathi--Laudenbach--Poenaru \cite{fathi_laudenbach_poenaru_79} or Farb--Margalit \cite{farb_margalit_11}).

\begin{definition}[Singular foliation]
A \emph{singular foliation} $F$ on a surface $(S, \ZZ)$ is a decomposition of $S$ into a disjoint union of subsets of $S$ (the \emph{leaves} of $F$) and a finite set of points of $S$ (the \emph{singular points} of $F$) such that the following holds:
\begin{enumerate}
\item For each nonsingular point $p \in S$, there is a smooth chart around $p$ that takes leaves to horizontal line segments. In addition, a transition map between any two charts is a smooth map of the form $(x, y) \mapsto (f(x,y), g(y))$.
\item For each singular point $p \in S$, there is a smooth chart around $p$ that takes leaves to the level sets of a $k$--pronged singularity with $k \geq 3$ or a one-pronged singularity (see \cref{fig:foliation}).
\end{enumerate}

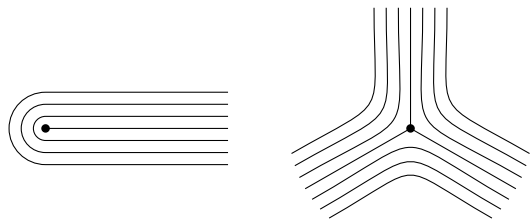
\begin{figure}[b]
 \centering
  \begin{tikzpicture}[scale=0.8]
   \fill (0,0) circle (2pt);
   \draw (0,0) -- (3,0);
   \draw (3,0.2) -- (0,0.2) arc (90:270:0.2) -- (3,-0.2);
   \draw (3,0.4) -- (0,0.4) arc (90:270:0.4) -- (3,-0.4);
   \draw (3,0.6) -- (0,0.6) arc (90:270:0.6) -- (3,-0.6);

   \begin{scope}[xshift=6cm]
    \fill (0,0) circle (2pt);
    \draw (0,0) -- (0,2);
    \draw (0,0) -- (-30:2);
    \draw (0,0) -- (-150:2);

    \draw (0.2,2) .. controls +(-90:2.3) and +(150:2.3) .. (-24:2);
    \draw (-0.2,2) .. controls +(-90:2.3) and +(30:2.3) .. (-156:2);
    \draw (-36:2) .. controls +(150:2.3) and +(30:2.3) .. (-144:2);

    \draw (0.4,2) .. controls +(-90:2.1) and +(150:2.1) .. (-18:2);
    \draw (-0.4,2) .. controls +(-90:2.1) and +(30:2.1) .. (-162:2);
    \draw (-42:2) .. controls +(150:2.1) and +(30:2.1) .. (-138:2);

    \draw (0.6,2) .. controls +(-90:1.9) and +(150:1.9) .. (-12:2);
    \draw (-0.6,2) .. controls +(-90:1.9) and +(30:1.9) .. (-168:2);
    \draw (-48:2) .. controls +(150:1.9) and +(30:1.9) .. (-132:2);
   \end{scope}
  \end{tikzpicture}
  \caption{Neighbourhoods of a one-pronged singularity and a three-pronged singularity of a measured foliation.}
  \label{fig:foliation}
\end{figure}
\end{definition}

\begin{definition}[Measured foliation]
 A \emph{measured foliation} on $S$ is a singular foliation $F$ together with a transverse measure $\mu$ such that every nonsingular point admits a smooth chart \hbox{$U \subset \mathbb R^2$} where $F$ is locally the pullback of the horizontal foliation on $\mathbb{R}^2$ and the measure $\mu$ is induced by~$|dy|$. Equivalently, if $\alpha, \beta: [0,1] \to S$ are two segments transverse to $F$ and isotopic through transverse segments whose endpoints always stay on the same leaf, then $\mu(\alpha) = \mu(\beta)$.
\end{definition}

\pagebreak

A measured foliation $F$ can be decomposed into a finite union of partial measured foliations $F'$ such that
\begin{itemize}
\item each $F'$ is supported on a connected subsurface $Y$ of $S$;
\item either every leaf of $F'$ is dense in $Y$, or $F'$ is a foliation of an annulus with all leaves closed and homotopic to the core curve.
\end{itemize}
 Each such partial measured foliation is called \emph{foliation component} of $F$. We assume that any annular components are maximal, that is, they contain all leaves of $F$ homotopic to its core curve. Thus, the supports of the foliation components are distinct and pairwise disjoint (up to~isotopy).

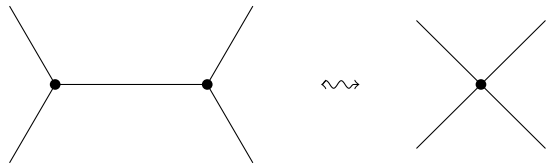
\begin{figure}
 \centering
  \begin{tikzpicture}
   \fill (0,0) circle (2pt);
   \draw (0,0) -- (2,0);
   \draw (0,0) -- (120:1.2);
   \draw (0,0) -- (-120:1.2);

   \fill (2,0) circle (2pt);
   \draw (2,0) -- ++(60:1.2);
   \draw (2,0) -- ++(-60:1.2);

   \draw[decorate, decoration={snake, segment length=2mm, amplitude=0.5mm}, <->] (3.5,0) -- (4,0);

   \begin{scope}[xshift=5.6cm]
    \fill (0,0) circle (2pt);
    \draw (0,0) -- (45:1.2);
    \draw (0,0) -- (135:1.2);
    \draw (0,0) -- (-135:1.2);
    \draw (0,0) -- (-45:1.2);
   \end{scope}
  \end{tikzpicture}
  \caption{A Whitehead move where a compact leaf is collapsed, combining two three-pronged singularities to a four-pronged singularity.}
  \label{fig:whitehead_isotopy}
\end{figure}

Two (measured) foliations on $(S, \ZZ)$ are \emph{Whitehead equivalent} if they differ by isotopy or Whitehead operations relative to $\ZZ$ (see \cref{fig:whitehead_isotopy}).
The set of Whitehead equivalence classes of measured foliations on the punctured surface $(S,\ZZ)$ is denoted by $\MF(S,\ZZ)$.
We also say that two measured foliations are \emph{topologically equivalent} if their underlying foliations are Whitehead equivalent.

\begin{definition}[Intersection number]
Let $(F,\mu)$ be a measured foliation and $\alpha$ be the homotopy class of a simple closed curve on $(S,\ZZ)$.
The \emph{intersection number of $\alpha$ with $(F,\mu)$} is defined as
$$i((F, \mu); \alpha) := \inf_{a \in \alpha} \int_a \mu.$$
\end{definition}

\begin{theorem}[Characterisation of Whitehead equivalence \cite{fathi_laudenbach_poenaru_79}]
Two measured foliations $(F, \mu)$ and $(F', \mu')$ on a surface~$(S,\ZZ)$ are Whitehead equivalent if and only if
$$i((F,\mu); \gamma) = i((F',\mu'); \gamma)$$
for every simple closed curve $\gamma$. \qed
\end{theorem}

Let $\cal S$ be the set of all simple closed curves. Then the intersection number is a map
$$ i \colon \MF(S, \ZZ) \times \mathcal{S} \to \R_{\geq 0},$$
which does not depend on the choice of the representative in the Whitehead equivalence class.
Equivalently, it gives a map $\MF(S, \ZZ) \hookrightarrow \R_{\geq 0}^{\cal S}$ from the space of measured foliation into the space of non-negative functions on $\cal S$. The map is injective and induces a natural topology on~$\MF(S, \ZZ)$.
For every simple closed curve $\alpha$, there exists a unique measured foliation $(F,\mu)\in \MF(S,\ZZ)$ such that $i(\alpha;\beta)=i((F,\mu);\beta)$ for each $\beta\in \mathcal S$.
Moreover, the intersection number extends to a continuous pairing
$$i: \MF(S,\ZZ) \times \MF(S,\ZZ) \to \mathbb R.$$
There is a natural action of $\mathbb R^+$  on $\MF(S, \ZZ)$ by multiplying the transverse measure by a constant.
The space of \emph{projective measured foliations} $\PMF(S, \ZZ) :=  (\MF(S, \ZZ) \setminus\{0\} ) /\mathbb{R}^+ $ is the quotient space under the action of scaling the transverse measure. The embedding  $\MF(S,\ZZ) \hookrightarrow \R^{\cal S}_{\geq 0}$ induces an embedding $\PMF(S,\ZZ) \hookrightarrow \mathbb {PR}^{\cal S}_{\geq 0}$ which shows that $\PMF(S, \ZZ)$ is compact.

\begin{definition}[Arational foliation]
 A measured foliation $(F, \mu)$ on $(S,\ZZ)$ is \emph{arational} if $i( (F, \mu); \gamma) > 0$ for every simple closed curve $\gamma$ on $(S,\ZZ)$. The space of \emph{arational (topological) foliations} $\mathcal{AF}(S,\ZZ)$ is obtained by taking the subspace of arational measured foliations in~$\MF(S,\ZZ)$ and forgetting the measure.
\end{definition}

If $(F, \mu)$ is arational then every half-leaf of $F$ that does not lead to a singularity is dense. In particular, the property of being arational only depends on the underlying topological foliation.
\begin{proposition}[Equivalence of arational foliations \cite{fathi_laudenbach_poenaru_79, rees_81}]\label{lem:arational_2}
Let $(F, \mu)$ and $(F', \mu')$ be measured foliations.
If $F$ is arational then it is topologically equivalent to $F'$ if and only if $i((F,\mu); (F',\mu')) = 0$. \qed
\end{proposition}

Given a half-translation surface $(S,q)$ and a slope $\theta\in\RP$, let $M_q(\theta)$ be the measured foliation on $(S,q)$ by straight lines with slope $\theta$,
where the transverse measure is locally the Euclidean distance between leaves.
Define $\MF(S,q;\ZZ) := M_q(\RP)$ to be the space of \emph{straight foliations} on $(S,q)$;
here, we are considering actual measured foliations, rather than their equivalence classes in $\MF(S,\ZZ)$.
Note that if $\theta \neq \theta'$ then $M_q(\theta)$ and $M_q(\theta')$ are transverse, and so their underlying topological foliations cannot be the same.

\begin{lemma}[Characterisation of arational foliations]\label{lem:arational}
The foliation $M_q(\theta)$ is arational on the topological surface $(S,\ZZ)$ if and only if there is no saddle connection on the half-translation surface~$(S,q)$ with slope $\theta$.
\end{lemma}

\proof
Suppose that $\alpha$ is a saddle connection on $(S,q)$ with slope $\theta$.
For $\epsilon > 0$, let $N_\epsilon(\alpha)$ be the~$\epsilon$--neighbourhood of $\alpha$ on $(S,q)$, and $\gamma_\epsilon$ be a component of $\partial N_\epsilon(\alpha)$. For $\epsilon$ sufficiently small,~$\gamma_\epsilon$ is an essential curve on $(S,\ZZ)$ that can be decomposed into at most two straight segments parallel to $\alpha$, together with at most two circular segments. Moreover, for $\epsilon > 0$ sufficiently small, there is a consistent choice of $\gamma_\epsilon$ such that all the $\gamma_\epsilon$ belong to a common homotopy class $\gamma$.
Since the straight segments are parallel to $\alpha$, the only contribution to $i((M_q(\theta), \mu); \gamma_\epsilon)$ comes from the circular segments. Each circular segment is a subpath of the boundary of the $\epsilon$--neighbourhood of a singularity, and so $\int_{\gamma_\epsilon} M_q(\theta) = 2k\epsilon$ for some (integer) constant $k$. Taking $\epsilon \rightarrow 0$, we see that $i((M_q(\theta), \mu); \gamma) = 0$ and so $M_q(\theta)$ is not arational.

For the converse, assume that $M_q(\theta)$ is not arational.
Then there exists a curve $\gamma$ on $(S,\ZZ)$ such that $i((M_q(\theta), \mu); \gamma) = 0$. By \cite[Proposition 5.7]{fathi_laudenbach_poenaru_79}, any geodesic representative of $\gamma$ has zero intersection number with $M_q(\theta)$.
A geodesic representative $\gamma_q$ is either a concatenation of saddle connections or a core curve of a Euclidean cylinder. In the first case, every saddle connection on~$\gamma_q$ has slope $\theta$, in the second case, every saddle connection on the boundary of the cylinder has slope $\theta$.
\endproof

Note that the assumption that we are in the fully-punctured setting is important here as a foliation of the form $M_q(\theta)$ might be arational on $S$ but not on $(S, \ZZ)$.
In the light of the lemma above, we call a slope $\theta \in \R P^1$ \emph{arational} if $M_q(\theta)$ contains no saddle connections.

\subsection{Combinatorial complexes}\label{subsec:combinatorial_complexes}

Here we recall various combinatorial complexes associated to surfaces that we use in this article.
The \emph{arc-and-curve graph}~$\AC(S, \ZZ)$ has as vertices the arcs or curves on $(S,\ZZ)$ (up to isotopy), with edges joining two arcs or curves if they have disjoint representatives. The \emph{arc graph} $\A(S,\ZZ)$ and \emph{curve graph} $\C(S, \ZZ)$ are the induced subgraphs of $\AC(S, \ZZ)$ whose vertices are respectively the arcs and the curves. We equip each of these graphs with the \emph{combinatorial metric}, where each edge is isometric to an interval of unit length.
Each of these graphs is connected, locally infinite, and has infinite diameter.

\begin{theorem}[Gromov hyperbolicity \cite{masur_minsky, masur_schleimer})]
The curve graph $\C(S, \ZZ)$ and the arc graph $\A(S, \ZZ)$ are Gromov hyperbolic. \qed
\end{theorem}

The inclusion $\A(S,\ZZ) \hookrightarrow \AC(S,\ZZ)$ is a quasi-isometric embedding if and only if $S$ has genus $0$ \cite{disarlo_0, schleimer}. In contrast, for the curve graph we have the following.

\begin{proposition}[Curve graph in the arc-and-curve graph \cite{kork_pap, schleimer}]\label{prop:Pap}
The natural inclusion \hbox{${\C(S,\ZZ) \hookrightarrow \AC(S,\ZZ)}$} is a quasi-isometry. Consequently, $\AC(S, \ZZ)$ is Gromov hyperbolic and there is a homeomorphism $\partial \C(S,\ZZ) \rightarrow \partial\AC(S, \ZZ)$ induced by the inclusion map. \qed
\end{proposition}

The boundary of the curve graph can be characterised explicitely in terms of arational foliations as follows.  An equivalent characterisation of $\partial\C(S,\ZZ)$ in terms of ending laminations is discussed by Hamenst\"adt \cite{hamenstaedt_06} and Pho-on \cite{pho-on_17}.

\begin{theorem}[Arational foliations form the boundary {\cite[Theorems 1.3 and 1.4]{klarreich_18}}]\label{thm:PO}
 There is a mapping class group equivariant homeomorphism $$K \from \partial\C(S,\ZZ)\to\AF(S, \ZZ).$$
Furthermore, if $\{\alpha_n\}$ is an admissible sequence in $\C(S, \ZZ)$
then, regarding the curves $\alpha_n$ as elements of $\PMF(S, \ZZ)$, every accumulation point of~$\{\alpha_n\}$ in~$\PMF(S, \ZZ)$ is topologically equivalent to the arational foliation $K([\{\alpha_n\}]) \in \AF(S, \ZZ)$. \qed
 \end{theorem}

The main object of interest in this paper is the \emph{saddle connection graph} $\A(S,q)$ associated to a half-translation surface $(S,q)$: this has the saddle connections on $(S,q)$ as vertices, with adjacency corresponding to pairs of non-crossing saddle connections. As we are working in the fully-punctured setting, the set of singularities of $(S,q)$ is equal to $\ZZ$, and so we may regard~$\A(S,q)$ as an induced subgraph of both $\A(S,\ZZ)$ and $\AC(S,\ZZ)$.

 Minsky and Taylor define the \emph{straightening map} $\straight : \AC(S, \ZZ) \to \A(S,q)$ as follows (they originally work with the arc graph, but their results also work in the setting of the arc-and-curve graph). Given an arc or curve $\alpha$, let $\straight(\alpha)$ be the set of all saddle connections that appear on the geodesic representative~$\alpha_q$. Observe that $\straight$ is the identity map when restricted to $\A(S,q)$, and so $\straight$ is a retraction.
 \cref{prop:MT} can now be reformulated in the language of combinatorial complexes. Recall that a \emph{clique} in a graph is a set of vertices that spans a complete subgraph.

\begin{proposition}[Straightening preserves disjointness \cite{minsky_taylor_17}] \label{prop:straightening_disjoint}
 If $\sigma \subset \AC(S,\ZZ)$ is a clique then $\straight(\sigma) \subset \A(S,q)$ is also a clique. \qed
\end{proposition}

Thus, $\straight : \AC(S, \ZZ) \to \A(S,q)$ is a coarse 1--Lipschitz retraction, in the sense that it sends sets of diameter at most 1 to sets of diameter at most 1. The following is immediate.

\begin{proposition}[Saddle connection graph in the arc-and-curve graph {\cite[Lemma 4.4]{minsky_taylor_17}}]\label{prop:MT_2}
The inclusion $\jmath \from \A(S,q)\to \AC(S,\ZZ)$ is an isometric embedding. Therefore, $\A(S,q)$ is Gromov hyperbolic and there exists a topological embedding
 $ \partial\jmath \from \partial\A(S,q)\to \partial\AC(S,\ZZ)~.$ \qed
\end{proposition}

\section{Paths in the saddle connection graph}

In this section, we find useful paths in the saddle connection graph. Our construction is inspired by the unicorn paths construction from \cite{hensel_przytycki_webb_15} in the topological arc graph $\A(S,\ZZ)$ which is based on Hatcher's surgery paths \cite{hatcher_91}.
We recall their construction here.

Let $\alpha$ and $\beta$ be arcs in minimal position on~$(S,\ZZ)$. Fix an orientation of $\alpha$ and $\beta$ and let $p$ and~$q$ respectively be the starting points of $\alpha$ and $\beta$.
Let $\alpha' \subset \alpha$ and $\beta' \subset \beta$ be subarcs starting at $p$ and $q$ respectively and ending at a common intersection point $z\in \alpha \cap \beta$.
If~$\alpha' \cup \beta'$ is an embedded arc then it is essential and we call it a \emph{unicorn arc obtained from~$\alpha$ and $\beta$}.
Unicorn arcs obtained from $\alpha$ and $\beta$ can be linearly ordered by declaring that $\alpha' \cup \beta' \leq \alpha'' \cup \beta''$ if and only if~$\beta' \subseteq \beta''$.
Let $(\gamma_1, . . . , \gamma_{n-1})$ be the ordered set of unicorn arcs.
A \emph{unicorn path} $P(\alpha, \beta)$ between~$\alpha$ and $\beta$ with starting points $p$ and $q$ is the sequence of vertices in $\A(S, \ZZ)$ given by
\[\alpha = \gamma_0, \gamma_1, . . . , \gamma_{n-1}, \gamma_n = \beta.\]
The homotopy class of $\gamma_i$ does not depend on the choice of representatives of $\alpha$ and~$\beta$.
Hensel--Przytycki--Webb showed that for any $\alpha, \beta \in \A(S, \ZZ)$ and any choice of orientations, the unicorn path $P(\alpha,\beta)$ is an (unparameterised) quasigeodesic in~$\A(S, \ZZ)$ \cite[Proposition~4.2]{hensel_przytycki_webb_15}.

\subsection{Bicorn paths}

We now generalise the construction of unicorn arcs to bicorn arcs.
Despite the similar name, this is a different construction to that of the bicorn \emph{curves} in \cite{przytycki_sisto_17}.

\begin{figure}[]
 \centering
  \begin{tikzpicture}[scale=0.8]
   \draw (0,0) -- (0,5);
   \draw[dotted] (0,5) -- (0,5.5);
   \draw (0,5.5) -- node[right]{$\alpha$} (0,6);

   \draw (-1,0.5) .. controls +(180:2) and +(180:3) ..
   (0,1.5) .. controls +(0:3) and +(0:3) ..
   (0,4) .. controls +(180:3) and +(180:3) ..
   (0,2.5) .. controls +(0:1) and +(225:0.5) .. (2,3);
   \draw (2.3,3.3) .. controls +(45:0.5) and +(270:0.5) .. (2.7,4.5);
   \draw[dotted] (2.7,4.5) -- node[right]{$\beta$} (2.7,5);

   \fill (0,0) circle (2pt) node[below]{$p_0$};
   \fill (0,1.5) circle (2pt) node[below left]{$p_1$};
   \fill (0,2.5) circle (2pt) node[below left]{$p_2$};
   \fill (0,4) circle (2pt) node[below left]{$p_3$};
   \fill (0,6) circle (2pt) node[above]{$p_{n+1}$};

   \fill (-1,0.5) circle (2pt);
  \end{tikzpicture}
  \caption{The arcs $\alpha$ and $\beta$ are in minimal position with intersection points $p_1,\ldots,p_n$.}
  \label{fig:setup_bicorn_arcs}
\end{figure}
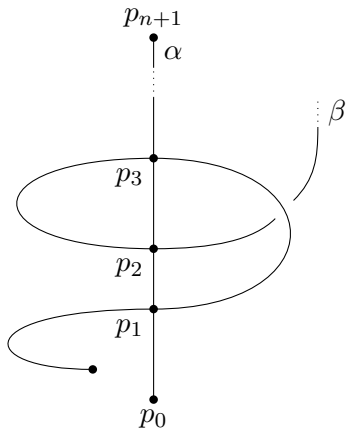

Let $\alpha$ and $\beta$ be arcs in minimal position on $S$.
Fix an orientation for $\alpha$, and let $p_0$ be its starting point.
Label the points in $\alpha \cap \beta$ by $p_1, \ldots, p_n$, where $n = i(\alpha, \beta)$,
so that the indices appear in increasing order along $\alpha$.
Let $p_{n+1}$ be the terminal endpoint of~$\alpha$ (see \cref{fig:setup_bicorn_arcs}).

For each $1 \leq i \leq n$, we construct the respective \emph{right} and \emph{left $(\alpha,\beta)$--bicorn arcs} $\gamma_i^+$ and~$\gamma_i^-$ as follows. Let $a_i \subseteq \alpha$ be the subsegment from $p_0$ to $p_i$.
Let $b_i^+ \subseteq \beta$ be the maximal subsegment of $\beta$ that begins at $p_i$,
contains no $p_j$ for $0 < j < i$ in its interior, and lies to the \emph{right} of $\alpha$ in a neighbourhood of $p_i$.
If the other endpoint of $b_i^+$ is an endpoint of $\beta$, then set $\gamma_i^+$ to be the $(\alpha, \beta)$--unicorn arc formed by concatenating $a_i$ with $b_i^+$. Otherwise, the other endpoint of~$\beta_i^+$ is~$p_j$ for some $0 < j < i$.
In this case, let $c_i^+ \subseteq a_i$ be the subsegment from $p_j$ to $p_0$,
and set~$\gamma_i^+$ to be the concatenation of $a_i$, $b_i^+$, and $c_i^+$.
The bicorn arc~$\gamma_i^-$ with its subsegments $b_i^-$ and~¸$c_i^-$ is defined similarly,
where we instead take~$b_i^-$ to lie to the \emph{left} of $\alpha$ near $p_i$.
Set $\gamma^\pm_0 = \beta$ and~$\gamma^\pm_{n+1} = \alpha$.
We shall respectively call $a_i$, $b^\pm_i$, and $c^\pm_i$ the $a$--, $b$--, and $c$--parts of $\gamma^\pm_i$.

\begin{lemma}[Ladder lemma] \label{lem:ladders}
 For each $0 \leq i \leq n$, the bicorn arcs $\gamma_i^+, \gamma_i^-, \gamma_{i+1}^+, \gamma_{i+1}^-$
 have pairwise disjoint representatives.
\end{lemma}

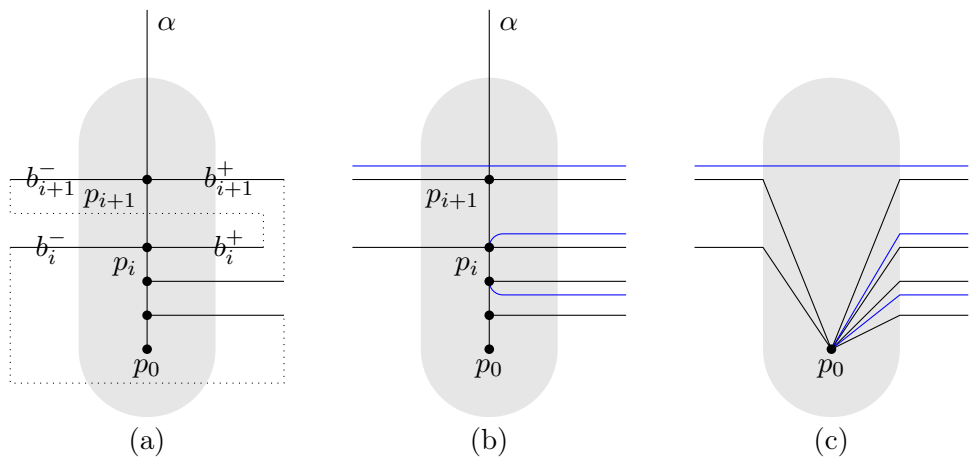
\begin{figure}
  \centering

  \begin{tikzpicture}[scale=0.9]
    \fill[gray!20](-1,0) arc (180:360:1cm)--
    (1,3) arc (0:180:1cm)--(-1,0);
    \draw (0,0)--(0,5)
    (0,0.5)--(2,0.5)
    (2,1)--(0,1)
    (-2,1.5)--(1.7,1.5)
    (-2,2.5)--(2,2.5);
    \draw[dotted](-2,2.5)--(-2,2.0)--
    (1.7,2.0)--(1.7,1.5)
    (2,2.5)--(2,1)
    (-2,1.5)--(-2,-0.5)--(2,-0.5)--(2,0.5);
    \draw (-1.4,2.5)node{$b^-_{i+1}$}
    (1.2,2.5)node{$b^+_{i+1}$}
    (-1.4,1.5)node{$b^-_{i}$}
    (1.2,1.5)node{$b^+_{i}$};
    \fill (0,0)circle(2pt) node[below]{$p_0$}
    (0,0.5)circle(2pt)
    (0,1.5)circle(2pt)node[below left]{$p_i$}
    (0,1)circle(2pt)
    (0,2.5)circle(2pt)node[below left]{$p_{i+1}$}
    (0,4.8)node[right]{$\alpha$};
     \draw (0,-1)node[below]{(a)};
    \begin{scope}[xshift=5cm]
        \fill[gray!20](-1,0) arc (180:360:1cm)--
    (1,3) arc (0:180:1cm)--(-1,0);
    \draw (0,0)--(0,5)
    (0,0.5)--(2,0.5)
    (2,1)--(0,1)
    (-2,1.5)--(2,1.5)
    (-2,2.5)--(2,2.5);
    \draw[blue](-2,2.7)--(2,2.7)
    (0,1) arc (-180:-90:0.2) -- (2,0.8)
    (0,1.5) arc (180:90:0.2) --(2,1.7);
    \fill (0,0)circle(2pt) node[below]{$p_0$}
    (0,0.5)circle(2pt)
    (0,1.5)circle(2pt)node[below left]{$p_i$}
    (0,1)circle(2pt)
    (0,2.5)circle(2pt)node[below left]{$p_{i+1}$}
    (0,4.8)node[right]{$\alpha$};
     \draw (0,-1)node[below]{(b)};
    \end{scope}
  \begin{scope}[xshift=10cm]
     \fill[gray!20](-1,0) arc (180:360:1cm)--
    (1,3) arc (0:180:1cm)--(-1,0);
    \draw (0,0)--(1,0.5)--(2,0.5)
     (0,0)--(-1,1.5)--(-2,1.5)
     (0,0)--(1,1)--(2,1)
     (0,0)--(1,1.5)--(2,1.5)
     (0,0)--(-1,2.5)--(-2,2.5)
    (0,0)--(1,2.5)--(2,2.5);
   \draw[blue](-2,2.7)--(2,2.7)
    (0,0)--(1,0.8)--(2,0.8)
    (0,0)--(1,1.7)--(2,1.7);
    \fill (0,0)circle(2pt) node[below ]{$p_0$};
    \draw (0,-1)node[below]{(c)};
  \end{scope}
  \end{tikzpicture}
  \caption{The shaded region in each figure represents $\mathscr R$. In (a), the $b$--part of $\gamma_i^+$ overlaps the $b$--parts of $\gamma^\pm_{i+1}$, i.e.\ $b^+_i=b^+_{i+1}\cup b^-_{i+1}$. In (b), we push the $b$--part of $\gamma^+_i$ to the blue arc so that it passes above $a_{i+1}$ near $p_{i+1}$. In (c), we perturb each connected component of the modified bicorn arcs inside $\mathscr R$.}\label{fig:ladderpath}
\end{figure}

\proof
We find disjoint representatives for all four arcs at once. First, we consider the $b$--parts of each bicorn arc. If the $b$--parts of some of the four arcs overlap, then $i<n$ and exactly one of $b^+_i$ or~$b^-_i$ equals $b^+_{i+1}\cup b^-_{i+1}$; without loss of generality suppose it is $b^+_i$. We homotope $\gamma^+_i$ by pushing its $b$--part off the unique side of $b^+_{i+1}\cup b^-_{i+1}$ so that it passes above $a_{i+1}$ near $p_{i+1}$ (see \cref{fig:ladderpath}\,(b)). Let $\hat{\gamma}^+_i$ be the resulting arc which is homotopic to $\gamma^+_i$, and $\hat{a}_{i+1}$ be the $a$--part of~$\hat{\gamma}^+_i$. If no two (original) $b$--parts share any overlap, we simply set $\hat{\gamma}^+_i=\gamma^+_i$.
In either case, the~$b$--parts of  $\hat{\gamma}^+_i,\gamma^-_i,\gamma^+_{i+1},\gamma^-_{i+1}$ have pairwise disjoint interiors.

Next, we consider the $a$--parts and $c$--parts. If $i=n$, let $\mathscr R$ be a small regular neighbourhood of $a_n$.
Otherwise, let $\mathscr R$ be a small regular neighbourhood of $\hat{a}_{i+1}$.
By choosing $\mathscr R$ small enough,
we may assume that $\mathscr R$ is simply connected,  and contains exactly one marked point $p_0$. Notice that  $\hat{\gamma}^+_i,\gamma^-_i$, $\gamma^+_{i+1},\gamma^-_{i+1}$ are disjoint outside $\mathscr R$. Inside $\mathscr R$, we homotope each  $\gamma\in\{\hat{\gamma}^+_i,\gamma^-_i, \gamma^+_{i+1},\gamma^-_{i+1}\}$ by perturbing every connected component of $(\gamma \setminus \partial \gamma)\cap \mathscr R$ as in \cref{fig:ladderpath}\,(c). The resulting four arcs are pairwise disjoint as desired.
\endproof

Consequently, we can define the respective \emph{left} and \emph{right bicorn paths}
\[\gamma_0^-, \ldots, \gamma_{n+1}^- \qquad\textrm{and}\qquad \gamma_0^+, \ldots, \gamma_{n+1}^+ \]
from $\beta$ to $\alpha$. Note that consecutive vertices in these paths may coincide.

\begin{lemma}[Bicorn arcs are near unicorn arcs]
 Each bicorn arc $\gamma_i^\pm$ has distance at most $1$ from some unicorn arc obtained from $\alpha$ and $\beta$ in $\A(S,\ZZ)$.
\end{lemma}

\proof
Note that $\gamma_0^\pm$ and $\gamma_1^\pm$ are themselves unicorn arcs. Now let $\gamma_i^\pm$ be a bicorn arc for~$i\geq 2$. Choose an orientation for $\beta$, and let $p_j$ be the first point that appears among $p_1, \ldots, p_{i}$ along~$\beta$.
The subsegment of $\beta$ from its starting point to~$p_j$ is either $b_j^+$ or $b_j^-$;
without loss of generality, suppose it is $b_j^+$.
Then $\gamma_j^+$ is a unicorn arc obtained from $\alpha$ and $\beta$.

If $j=i$, i.e.\ $\gamma_i^+$ is itself a unicorn arc, the statement follows from \cref{lem:ladders}.
If $j<i$, then the (open) $b$--parts of $\gamma_i^\pm$ and $\gamma_j^+$ do not overlap each other. By perturbing each of $\gamma^+_j$ and $\gamma^\pm_i$ inside a small regular neighbourhood of $a_i$ as in the proof of \cref{lem:ladders}, we see that $i(\gamma_i^\pm, \gamma_j^+) = 0$.
\endproof

Combining the previous lemma with \cite[Proposition 4.2]{hensel_przytycki_webb_15}, we deduce the following.

\begin{corollary}[Bicorn paths are quasigeodesics]
 For all $\alpha, \beta \in \A(S,\ZZ)$, the left and right bicorn paths from $\beta$ to $\alpha$ are (unparameterised) quasigeodesics in $\A(S,\ZZ)$. \qed
\end{corollary}

\subsection{From bicorn paths in \texorpdfstring{$\A(S,\ZZ)$}{A(S,Z)} to ladder paths in \texorpdfstring{$\A(S,q)$}{A(S,q)}}

We now exploit the construction of bicorn paths in the arc graph to obtain paths with similar properties in the saddle connection graph.

Consider two saddle connections $\alpha$ and $\beta$. For convenience in the proofs, we assume that they are vertical and horizontal, respectively.
Let $\gamma_0^-, \ldots, \gamma_{n+1}^-$ and $\gamma_0^+, \ldots, \gamma_{n+1}^+$ be the two bicorn paths from $\beta$ to $\alpha$ as in the previous section.
The geodesic representative of each (oriented) bicorn arc $\gamma_i^\pm$ is a concatenation of saddle connections; let $\delta_i^\pm$ be the first saddle connection (starting from $p_0$) appearing on the geodesic representative.
We prove that the sequences $\delta_i^\pm$ inherit several properties from the bicorn arcs~$\gamma_i^\pm$, as well as enhancing them.

\begin{proposition}[Properties of ladders paths] \label{lem:straight_ladders}
 Let $\alpha, \beta \in \A(S,q)$ be respectively vertical and horizontal, and set $n = i(\alpha, \beta)$.
 Then the sequences $\delta_0^-, \ldots, \delta_{n+1}^-$ and $\delta_0^+, \ldots, \delta_{n+1}^+$ defined as above form a pair of paths (where consecutive vertices may coincide) in $\A(S,q)$ with the following properties:
 \begin{itemize}
  \item $\delta_0^\pm = \beta$ and $\delta_{n+1}^\pm = \alpha$,
  \item $\delta_i^\pm$ has no transverse intersections with any of $\delta_i^\mp$, $\delta_{i+1}^+$, and $\delta_{i+1}^-$ for all $0 \leq i \leq n$,
  \item $\delta_i^+$ has positive slope and $\delta_i^-$ has negative slope for all $1 \leq i \leq n$, and
  \item the slopes of $\delta_0^+,\ldots,\delta_{n+1}^+$ are monotonically increasing and the slopes of $\delta_0^-,\ldots,\delta_{n+1}^-$ are monotonically decreasing.
 \end{itemize}
\end{proposition}

\proof
The first property follows immediately from the definition of the $\delta_i^\pm$.
The second property is a consequence of \cref{lem:ladders} and \cref{prop:straightening_disjoint}. The same argument also shows that $\delta_0^-, \ldots, \delta_{n+1}^-$ and $\delta_0^+, \ldots, \delta_{n+1}^+$ form paths (possibly where consecutive vertices coincide).
For the third and fourth properties, the proofs for the $\delta_i^+$ and $\delta_i^-$ are similar. Without loss of generality, we restrict attention to the sequence $\delta_i^+$ (defined using the right bicorn path).

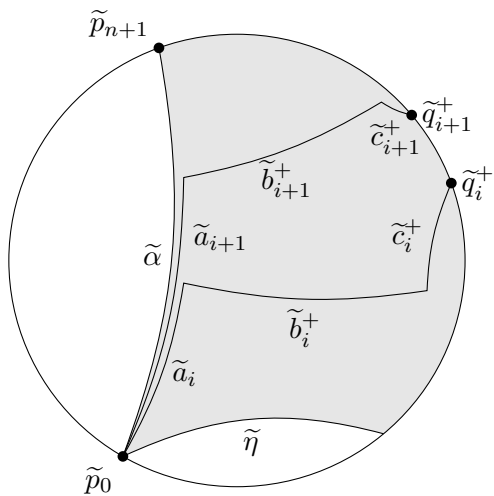
\begin{figure}
 \centering
  \begin{tikzpicture}
   \fill[color=gray!20] (-120:3)
   to[bend left=20] (-50:3)
   arc (-50:110:3)
   to[bend left=15] (-120:3);

   \draw (0,0) circle (3cm);

   \draw (-120:3) to[bend right=15] node[left]{$\widetilde{\alpha}$} (110:3);

   \fill (-120:3) circle (2pt) node[below left]{$\widetilde{p}_0$};
   \fill (110:3) circle (2pt) node[above left]{$\widetilde{p}_{n+1}$};

   \draw (-120:3) to[bend left=20] node[below]{$\widetilde{\eta}$} (-50:3);

   \draw (-120:3) to[bend right=10] node[right]{$\widetilde{a}_i$}
   (-0.7,-0.3) to[bend right=10] node[below]{$\widetilde{b}_i^+$}
   (2.5,-0.4) to[bend left=10] node[left]{$\widetilde{c}_i^+$}
   (20:3);
   \fill (20:3) circle (2pt) node[right]{$\widetilde{q}_i^+$};

   \draw (-120:3) to[bend right=11]
   (-0.7,1.1) to[bend right=10] node[below]{$\widetilde{b}_{i+1}^+$}
   (1.9,2.1) to[bend right=10] node[below]{$\widetilde{c}_{i+1}^+$}
   (40:3);
   \draw (-0.7,0.3) node[right]{$\widetilde{a}_{i+1}$};
   \fill (40:3) circle (2pt) node[right]{$\widetilde{q}_{i+1}^+$};
  \end{tikzpicture}
  \caption{Lifting the bicorn arcs $\gamma_i^+$ and $\gamma_{i+1}^+$ to $\mathbf{X}$; with the corresponding segments $\widetilde{a}_i$, $\widetilde{b}_i^+$,\ldots\ 
  All corners in $\hat X$ appearing in the figure are right angles.}
  \label{fig:ladder_paths_monotonicity}
\end{figure}

 Choose a lift $\widetilde{\alpha}$ of $\alpha$ to the Minsky--Taylor compactification $\mathbf{X}$ as described in \cref{thm:compactification} (see \cref{fig:ladder_paths_monotonicity} for the setup).
 Let $\widetilde{p}_0$ and $\widetilde{p}_{n+1}$ be the corresponding lifts of the endpoints~$p_0$ and~$p_{n+1}$.
 Let $\widetilde{\gamma}_i^+$ be the lift of $\gamma_i^+$ that shares a common initial segment with $\widetilde{\alpha}$ and let $l_i^+$ be its (oriented) geodesic representative. Note that the first saddle connection appearing on $l_i^+$ is a lift of $\delta_i^+$; we shall denote this by $\widetilde{\delta}_i^+$.
 Next, let $\widetilde{\eta}$ be the (shortest) complete geodesic starting at $\widetilde{p}_0$ and going to the right such that the angle at $\widetilde{p}_0$ between $\widetilde{\eta}$ and~$\widetilde{\alpha}$ is~$\sfrac{\pi}{2}$.
 Let $\mathfrak{C}=H_R(\widetilde{\alpha})\cap H_L(\widetilde{\eta})$ be the intersection of the \emph{right} side of $\widetilde{\alpha}$ and the \emph{left} side of $\widetilde{\eta}$ (as defined in \cref{ssec:MinskyTaylor}).
 By \cref{lem:convex}, both~$H_R(\widetilde{\alpha})$ and $H_L(\widetilde{\eta})$ are convex, hence $\mathfrak{C}$ is also convex.
 Observe that a small metric neighbourhood of $\widetilde{p}_0$ in $\mathfrak{C}$ is isometric to a Euclidean sector of angle $\frac{\pi}{2}$.
 By the construction of bicorn paths, $\widetilde{\gamma}_i^+$ is contained in~${\mathfrak{C}}$ for each $1\leq i\leq n+1$. As~$\mathfrak{C}$ is convex, it follows that $l_i^+ \subset \mathfrak{C}$.
 In particular, the slope of $\widetilde{\delta}_i^+$, hence $\delta_i^+$, is positive for $1 \leq i \leq n$ which proves the third property.

 Now, let $\widetilde{q}_i^+$ be the terminal endpoint of $\widetilde{\gamma}_i^+$. As the arcs~$\widetilde{\gamma}_i^+$ and~$\widetilde{\gamma}_{i+1}^+$ have no transverse intersections, the points~ $\widetilde{p}_0, \widetilde{q}_i^+, \widetilde{q}_{i+1}^+, \widetilde{q}_{n+1}^+=\widetilde{p}_{n+1}$ appear in the given (anticlockwise) cyclic order on the boundary $\partial\mathbf{X}\simeq {S}^1$ (see \cref{thm:compactification}).
 Therefore $\tilde{p}_0$ and $\tilde{q}_{i+1}^+$ both lie in $H_R(\widetilde{\alpha})\cap H_L(l_i^+)$.
 Again by convexity, and the fact that $l_{i+1}^+$ is the geodesic from $\tilde{p}_0$ to $\tilde{q}_{i+1}^+$, we have
 \[\widetilde{\delta}_{i+1}^+ \subseteq l_{i+1}^+ \subset H_R(\widetilde{\alpha})\cap H_L(l_i^+) \subseteq \mathfrak{C}. \]
  Since $\widetilde{\delta}_i^+$ is the first saddle connection appearing on $l_i^+$, it follows that the slope of $\widetilde{\delta}_{i+1}^+$ is at least that of $\widetilde{\delta}_i^+$, yielding the fourth property as desired.
\endproof

We call the paths
\[\beta = \delta_0^-, \ldots, \delta_{n+1}^- = \alpha \quad \textrm{and} \quad \beta = \delta_0^+, \ldots, \delta_{n+1}^+ = \alpha\]
respectively the \emph{left} and \emph{right ladder paths} from $\beta$ to $\alpha$ in $\A(S,q)$. Note that consecutive vertices can coincide, i.e.\ we can have $\delta_i^\pm = \delta_{i+1}^\pm$; this is illustrated in the following example.

\begin{example}[Ladder paths for the once-marked torus]
 For the once-marked torus, the arc graph and the saddle connection graph are both isomorphic to the Farey tesselation graph (see \cite[Example 3.6]{disarlo_randecker_tang_18}).
 We shall view vertices as slopes in $\RP$ using the standard picture of the Farey tesselation.
 Using this, we can give an explicit description of ladder paths between two saddle connections~$\alpha$ and $\beta$.

 If~$\alpha$ and~$\beta$ are adjacent, then the ladder path is the edge connecting them, so we shall assume otherwise. Consider the set of all triangles in the Farey tesselation that contain an edge separating~$\alpha$ and $\beta$ on $\RP$. The union of these triangles is a polygon $P$ with~$\alpha$ and~$\beta$ among its vertices.
 Observe that $\{\alpha, \beta\}$ divides $\RP$ into two subintervals~$I^+$ and~$I^-$ such that all $\delta^\pm_i$ lie in~$I^\pm$. For each $i$, the vertices~$\delta_i^+$ and~$\delta^-_i$ are adjacent, and so all vertices along the ladder paths must appear as vertices on $P$.
 On the other hand, there is a unique path from $\beta$ to $\alpha$ such that all intermediate vertices appear on $P$ and lie in $I^\pm$; this path must be one of the two paths along the boundary of $P$ from $\beta$ to $\alpha$.
 Therefore, the ladder paths are the two paths along the boundary of $P$ that connect $\beta$ to~$\alpha$ (allowing for consecutive vertices to coincide).
 An explicit example is given in \cref{fig:exa_farey_tesselation}.

\begin{figure}
 \centering
  \begin{tikzpicture}[scale=0.65]
   \draw[thick, fill=gray!20] (-118:5)
   to[out=62, in=90, looseness=1.2] (-90:5)
   arc (180:90:5)
   to[out=180, in=217, looseness=1.3] (37:5)
   to[out=217, in=-90, looseness=1.2] (90:5)
   arc (0:-90:5)
   to[out=0, in=37, looseness=1.3] (-143:5)
   (-143:5) to[out=37, in=53, looseness=1.5] (-127:5)
   to[out=53, in=62, looseness=1.5] (-118:5);

   \draw[densely dashed, thick, color=gray] (-118:5) to[out=70, in=-143] (37:5);

   \draw (0,0) circle (5cm);
   \fill (0:5) circle (2pt) node[right]{$1$};
   \fill (21:5) circle (2pt) node[right]{$\frac{3}{2}$};
   \fill (37:5) circle (3pt) node[right]{$2$};
   \fill (53:5) circle (2pt) node[above right]{$3$};

   \fill (90:5) circle (2pt) node[above]{$\frac{1}{0}$};

   \fill (159:5) circle (2pt) node[left]{$-\frac{3}{2}$};
   \fill (143:5) circle (2pt) node[left]{$-2$};
   \fill (127:5) circle (2pt) node[above left]{$-3$};

   \fill (180:5) circle (2pt) node[left]{$-1$};

   \fill (-159:5) circle (2pt) node[left]{$-\frac{2}{3}$};
   \fill (-143:5) circle (2pt) node[left]{$-\frac{1}{2}$};
   \fill (-127:5) circle (2pt) node[below left]{$-\frac{1}{3}$};
   \fill (-118:5) circle (3pt) node[below]{$-\frac{1}{4}$};

   \fill (-21:5) circle (2pt) node[right]{$\frac{2}{3}$};
   \fill (-37:5) circle (2pt) node[right]{$\frac{1}{2}$};
   \fill (-53:5) circle (2pt) node[below right]{$\frac{1}{3}$};
   \fill (-90:5) circle (2pt) node[below]{$0$};

   \draw (-90:5) -- (90:5);
   \draw (0:5) arc (-90:-180:5)
    arc (0:-90:5)
    arc (90:0:5)
    arc (180:90:5);

   \draw (0:5) to[out=180, in=217, looseness=1.3] (37:5)
   to[out=217, in=-90, looseness=1.2] (90:5);
   \draw (0:5) to[out=180, in=143, looseness=1.3] (-37:5)
   to[out=143, in=90, looseness=1.2] (-90:5);
   \draw (180:5) to[out=0, in=-37, looseness=1.3] (143:5)
   to[out=-37, in=-90, looseness=1.2] (90:5);
   \draw (180:5) to[out=0, in=37, looseness=1.3] (-143:5)
   to[out=37, in=90, looseness=1.2] (-90:5);

   \draw (37:5) to[out=217, in=-159, looseness=1.5] (21:5)
   to[out=-159, in=180, looseness=1.4] (0:5);
   \draw (-37:5) to[out=143, in=159, looseness=1.5] (-21:5)
   to[out=159, in=180, looseness=1.4] (0:5);
   \draw (143:5) to[out=-37, in=-21, looseness=1.5] (159:5)
   to[out=-21, in=0, looseness=1.4] (180:5);
   \draw (-143:5) to[out=37, in=21, looseness=1.5] (-159:5)
   to[out=21, in=0, looseness=1.4] (180:5);

   \draw (37:5) to[out=217, in=-127, looseness=1.5] (53:5)
   to[out=-127, in=-90, looseness=1.3] (90:5);
   \draw (-37:5) to[out=143, in=127, looseness=1.5] (-53:5)
   to[out=127, in=90, looseness=1.3] (-90:5);
   \draw (143:5) to[out=-37, in=-53, looseness=1.5] (127:5)
   to[out=-53, in=-90, looseness=1.3] (90:5);
   \draw (-143:5) to[out=37, in=53, looseness=1.5] (-127:5)
   to[out=53, in=90, looseness=1.3] (-90:5);

   \draw (-127:5) to[out=53, in=62, looseness=1.5] (-118:5)
   to[out=62, in=90, looseness=1.2] (-90:5);
  \end{tikzpicture}
  \caption{The ladder paths from the saddle connection with slope $-\frac{1}{4}$ to the one with slope $2$ are the two paths that bound the gray polygon.}
  \label{fig:exa_farey_tesselation}
\end{figure}
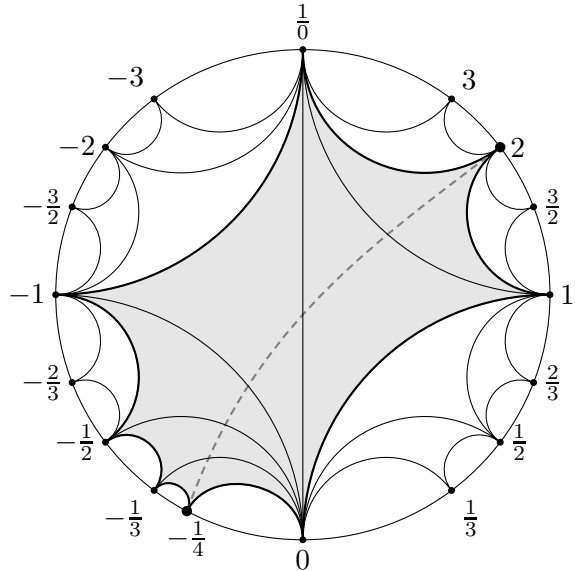

 Note that for the once-marked torus, the ladder paths in the saddle connection graph coincide with the bicorn paths in the arc graph.
\end{example}

In the special case of the Farey tesselation graph, our ladder paths recover the ``ladders'' from \cite{baik_kim_kwak_shin_19}, also called ``fans'' in \cite{hatcher_17}.

\pagebreak

\subsection{Ladder paths are close to every path}

This subsection is dedicated to proving the following proposition.

\begin{proposition}[Ladder paths are bottlenecks]\label{prop:4-neighbourhood}\label{prop:ladder_paths_close_to_every_path}
  Let $\alpha$ and $\beta$ be distinct saddle connections and suppose $\beta = \gamma_0, \gamma_1, \ldots, \gamma_m = \alpha$ is an arbitrary path in $\A(S,q)$ from $\beta$ to $\alpha$.
  Then the ladder paths constructed in \cref{lem:straight_ladders} are both contained in the $3$--neighbourhood of the given path.
\end{proposition}

This proposition would be enough to prove that $\A(S,q)$ is a quasitree using Manning's Bottleneck Criterion \cite{manning_05}. However, since our goal is to show that $\A(S,q)$ is quasi-isometric to the specific tree $T_\infty$, we shall construct the quasi-isometry explicitly in the next section.

To prove \cref{prop:ladder_paths_close_to_every_path}, we require a key lemma which is also useful for proving \cref{thm:main}.

\begin{lemma}[Linking slopes]\label{lem:separating_slopes}
 Let $\alpha_1, \alpha_2, \alpha_3, \alpha_4 \in \A(S,q)$ be saddle connections
 whose slopes are distinct and appear in the given cyclic order in $\RP$.
 Assume that $\alpha_1$ and $\alpha_3$ do not cross, and that $\alpha_2$ and $\alpha_4$ do not cross.
 Then $d_{\A(S,q)}(\{\alpha_1,\alpha_3\},\{\alpha_2,\alpha_4\})\leq 2$.
\end{lemma}
 \proof
  By applying three $\SL$--deformations in turn, we may assume that $\alpha_1$ is vertical (by rotation), $\alpha_3$ is horizontal (by horocycle flow), and that $\alpha_2$ and $\alpha_4$ are perpendicular on~$(S,q)$ (by Teichmüller flow).
  Let $\beta$ be a shortest saddle connection on $(S,q)$ and $p$ be its midpoint. Let~$r$ be the length of $\beta$ and~$U$ be the $\frac{r}{2}$--neighbourhood of~$p$ on $(S,q)$. Since $\beta$ is a shortest saddle connection,~$U$ is isometric to a Euclidean disc. Since $\alpha_1$ and $\alpha_3$ are perpendicular, they cannot both cross $\beta$, for otherwise
  $\alpha_1$ and $\alpha_3$ would intersect inside $U$. Therefore, $d_{\A(S,q)}(\beta,\{\alpha_1,\alpha_3\})\leq 1$. Similarly, we have $d_{\A(S,q)}(\beta,\{\alpha_2,\alpha_4\})\leq 1$. Hence $d_{\A(S,q)}(\{\alpha_1,\alpha_3\},\{\alpha_2,\alpha_4\})\leq 2$.
 \endproof

\begin{proof}[Proof of \cref{prop:4-neighbourhood}]
 If $d(\alpha, \beta) = 1$, the ladder paths both coincide with the edge between~$\alpha$ and~$\beta$. Thus, each ladder path is in the $1$--neighbourhood of $\{\alpha, \beta\}$, hence in the $1$--neighbourhood of any path from $\alpha$ to $\beta$.

 If $d(\alpha, \beta) > 1$ then $\alpha$ and $\beta$ cannot be parallel. Therefore, we can apply an $\SL$--deformation to arrange for $\alpha$ and $\beta$ to be vertical and horizontal respectively. Let $\beta = \delta_0^-, \ldots, \delta_{n+1}^- = \alpha$ and $\beta = \delta_0^+, \ldots, \delta_{n+1}^+ = \alpha$ be the ladder paths as constructed in \cref{lem:straight_ladders}.
 We shall prove that for every $i\in \{1,\ldots,n\}$, there exists some saddle connection $\gamma_j$ along the given path that lies within distance~$3$ of both $\delta_i^+$ and $\delta_i^-$.
 If some $\gamma_j$ is parallel to $\delta_i^+$ or $\delta_i^-$ then we are done, so we may suppose otherwise.

 Observe that the slopes of $\delta_i^+$ and $\delta_i^-$ cut $\RP$ into exactly two components, one containing the slope of $\beta=\gamma_0$ and the other containing the slope of $\alpha=\gamma_m$. Since $\gamma_0, \gamma_1, \ldots, \gamma_m$ is a path from $\beta$ to $\alpha$, there exists some $0 \leq j < m$ such that the slopes of $\gamma_j$ and $\gamma_{j+1}$ lie in different components.
 This means that the slopes of $\gamma_j$, $\delta_i^+$, $\gamma_{j+1}$, and $\delta_i^-$ appear in the given cyclic order along $\RP$. Applying \cref{lem:separating_slopes}, we deduce that $d_{\A(S,q)}(\{\gamma_j,\gamma_{j+1}\},\{\delta_i^+,\delta_i^-\})\leq 2$. As $d_{\A(S,q)}(\delta_i^+,\delta_i^-) = 1$, we have that $\delta_i^+$ and to~$\delta_i^-$ both lie within distance $3$ of some vertex along the given path.
\end{proof}

\cref{prop:ladder_paths_close_to_every_path} can also be used to deduce Gromov hyperbolicity of $\A(S,q)$. This would also follow later from \cref{thm:main}, however, we include a proof here to obtain explicit constants. Our proof is essentially the same as the proof of uniform hyperbolicity of the arc graph in \cite{hensel_przytycki_webb_15}.

\begin{proposition}[$4$--hyperbolicity of the saddle connection graph] \label{prop:hyperbolicity_constant_saddle_connection_complex}
 The saddle connection \linebreak graph $\A(S,q)$ is $4$--hyperbolic in the sense that every geodesic triangle has a $4$--centre.

\begin{proof}
 Let $\alpha_1$, $\alpha_2$, and $\alpha_3$ be a triple of saddle connections on $(S,q)$.
 By \cite[Lemma 3.4]{hensel_przytycki_webb_15}, for each $i=1,2,3$, there is a unicorn arc $\gamma_i$ obtained from $\alpha_i$ and $\alpha_{i+1}$ (with indices considered modulo $3$), such that $\gamma_1,\gamma_2$ and $\gamma_3$ are pairwise disjoint.
 As every unicorn arc is also a bicorn arc,~$\gamma_i$ is contained in some bicorn path between $\alpha_i$ and $\alpha_{i+1}$
 in $\A(S,\ZZ)$.
 Using \cref{lem:straight_ladders}, from the bicorn path, we get a ladder path in the saddle connection graph $\A(S,q)$ connecting $\alpha_i$ and $\alpha_{i+1}$. Let $\delta_i$ be the saddle connection along the ladder path
 corresponding to the unicorn arc~$\gamma_i$. Then $\delta_1$, $\delta_2$, and $\delta_3$ are pairwise disjoint by \cref{prop:straightening_disjoint}.

 Now let $[\alpha_1,\alpha_2]\cup[\alpha_2,\alpha_3]\cup[\alpha_3,\alpha_1]$ be an arbitrary geodesic triangle in $\A(S,q)$ with vertices $\alpha_1,\alpha_2$ and $\alpha_3$. Since $\{\delta_1, \delta_2, \delta_3\}$ has diameter at most 1 in $\A(S,q)$, it follows using \cref{prop:4-neighbourhood} that $\delta_1$ lies within distance 4 of each side $[\alpha_1,\alpha_2]$, $[\alpha_2,\alpha_3]$, and $[\alpha_3,\alpha_1]$.
\end{proof}
\end{proposition}

\section{Quasi-isometry via the graph of slopes}

In this section, we construct an explicit quasi-isometry from the saddle connection graph to the regular infinite-valent tree. This map is defined by a quotient map that factors through the graph of slopes.

\subsection{Graph of slopes}

Let us introduce the main object of interest in this section.

\begin{definition}[Graph of slopes]
 The \emph{graph of slopes} $\G(S,q)$ associated to a half-translation surface $(S,q)$ has
 as vertices slopes $\theta \in \RP$ admitting a saddle connection on $(S,q)$,
 with two slopes $\theta \neq \theta'$ declared adjacent if they can be realised by a pair of disjoint saddle connections.
\end{definition}

There is a natural map $\Theta \from \A(S,q) \to \G(S,q)$ sending a saddle connection to its slope;
this is the quotient map defined on $\A(S,q)$ where the equivalence classes are (maximal) sets of parallel saddle connections.
We use the following observation to show that this map is a quasi-isometry.

\begin{lemma}[Quotient maps with bounded equivalence classes are quasi-isometries] \label{lem:quotient_maps_quasi_isometries}
 Let $G$ be a graph and $\sim$ be an equivalence relation on the vertex set of $G$.
 If each equivalence class has diameter at most $K \geq 0$ then the quotient map $\pi: G \to G/\sim$ is a $(K+1,K)$--quasi-isometry. \qed
\end{lemma}

\begin{corollary}[Slope map is quasi-isometry]
 The slope map $\Theta \from \A(S,q) \to \G(S,q)$ is a $(2,1)$--quasi-isometry.

\begin{proof}
 Each equivalence class $\Theta^{-1}(\theta)$ is a set of parallel saddle connections,
 and so has diameter at most $1$ in $\A(S,q)$.
\end{proof}
\end{corollary}

We can use the slope map to define ladder paths in the graph of slopes, similarly to ladder graphs in the saddle connection graph:
For two slopes $\theta, \theta' \in \G(S,q)$, we define (pairs of) \emph{ladder paths in the graph of slopes} from $\theta$ to~$\theta'$ to be images under $\Theta$ of pairs of right and left ladder paths from any $\alpha \in \Theta^{-1}(\theta)$ to any~$\alpha' \in \Theta^{-1}(\theta')$ in~$\A(S,q)$.
The slopes along each of these ladder paths in the graph of slopes appear in cyclic order on~$\RP$; anticlockwise for the right ladder path, clockwise for the left ladder path.

Since the slope map is $1$--Lipschitz, we immediately deduce a version of \cref{prop:ladder_paths_close_to_every_path} for the graph of slopes $\G(S,q)$.

\begin{corollary}[Bottlenecks in the graph of slopes] \label{prop:slopes_ladder_paths_close_to_every_path}
 Let $\theta, \theta'$ be two slopes with $d(\theta,\theta') \geq 1$ in~$\G(S,q)$ and let $$\theta = \theta_0, \theta_1, \ldots, \theta_m = \theta'$$ be a path in $\G(S,q)$ connecting $\theta$ and $\theta'$.
 Then the ladder paths in the graph of slopes are both contained in the $3$--neighbourhood of the given path. \qed
\end{corollary}

\subsection{Balls in the graph of slopes}

 We now investigate the behaviour of balls in $\G(S,q)$ in conjunction with the topology of~$\RP$.
 A~set of vertices $A \subseteq \G(S,q)$ is called
 \emph{closed} if it is closed as a subset of $\RP$.
 Given $\theta \in \G(S,q)$,~let
  \[B(\theta, r) := \set{\theta' \in \G(S,q) \st d_{\G(S,q)}(\theta, \theta') \leq r}\]
 be the ball of radius $r\geq 0$.
 The following propositions are important for our proof of \cref{thm:main}.

 \begin{proposition}[Balls are closed] \label{prop:balls_in_slope_graph_are_closed}
  For any $\theta\in \G(S,q)$ and $r \geq 0$, the ball $B(\theta, r)$ is closed in~$\RP$.
 \end{proposition}

 Since balls in $\G(S,q)$ determine closed and countable subsets of $\RP$, it follows that
 $\RP \setminus B(\theta, r)$ has full measure and is a disjoint union of open intervals.
 We show that each such interval is infinitely subdivided when the radius of the ball is increased by $1$.

 \begin{proposition}[Infinite subdivision] \label{prop:infinitely_many_r+1}
 Given $\theta\in \G(S,q)$ and $r \geq 0$, every connected component of $\RP \setminus B(\theta, r)$ is an open interval containing infinitely many slopes in~$B(\theta, r+1)$.
 \end{proposition}

 To prove the above, we require several convergence properties for sequences of saddle connections and slopes.
 The following is from \cite[Lemma 5.6]{disarlo_randecker_tang_18}.
 Recall that~$M_q(\theta)$ is the measured foliation on $(S,q)$ by straight lines of slope~$\theta \in \RP$, with the measure given locally by Euclidean distance between leaves.

 \begin{lemma}[Saddle connections with convergent slopes \cite{disarlo_randecker_tang_18}]
  Let $\alpha_n \in \A(S,q)$ be a sequence of saddle connections whose slopes $\theta_n$ converge to $\theta\in\RP$.
  Then, upon passing to a subsequence, there exist lifts $\widetilde{\alpha}_n$ to the universal cover $(\widetilde{S}, \widetilde{q})$ converging to a geodesic ray $\widetilde{l}$ in the following sense: for every metric ball $\widetilde{B} \subset (\widetilde{S}, \widetilde{q})$, the sequence $\widetilde{\alpha}_n \cap \widetilde{B}$ converges to $\widetilde{l} \cap \widetilde{B}$ with respect to the Hausdorff topology on $\widetilde{B}$.
  Furthermore, $\widetilde{l}$ descends to a concatenation of separatrices $l$ of $M_q(\theta)$ on $(S,q)$. \qed
 \end{lemma}

 Note that for the above lemma, we are taking the universal cover of $(S,q)$ without puncturing at the singularities.
 Thus, any sufficiently small (closed) metric ball is topologically a closed disc.

 \begin{lemma}[Convergent slopes and disjointness] \label{lem:convergence_and_disjointness}
  Suppose that $\alpha_n, \beta_n \in \A(S,q)$ are sequences of saddle connections whose respective slopes
  $\theta_n$, $\theta'_n$ converge to distinct slopes $\theta$, $\theta' \in \RP$. Assume that $\alpha_n$ and $\beta_n$ are disjoint for all~$n$.
  Then there exists a pair of disjoint saddle connections $\alpha$ and $\beta$ with slopes $\theta$ and $\theta'$.
 \end{lemma}

 \proof
 By passing to appropriate subsequences, there exist leaves $l, l'$ of $M_q(\theta)$ and $M_q(\theta')$ respectively satisfying the above lemma for the sequences $\alpha_n$ and $\beta_n$.

 We claim that $l$ and~$l'$ do not intersect transversely. If they do, then, up to deck transformations, the associated geodesic rays $\widetilde{l}$ and $\widetilde{l'}$ intersect transversely at some point $\widetilde{p}$ on the universal cover. Take $\widetilde{B}$ to be a (closed) metric ball centred at $\widetilde{p}$, chosen sufficiently small so that it is a closed disc containing no endpoint of $\widetilde{l}$ nor $\widetilde{l'}$.
 Since $\widetilde{l} \cap \widetilde{B}$ and $\widetilde{l'} \cap \widetilde{B}$ intersect transversely, $\widetilde{l} \cap \partial \widetilde{B}$ cuts $\partial \widetilde{B}$ into two connected components such that the two points in~$\widetilde{l'} \cap \partial\widetilde{B}$ lie in different components. For appropriate lifts $\widetilde{\alpha}_n$ and $\widetilde{\beta}_n$, the sets $\widetilde{\alpha}_n \cap \partial \widetilde{B}$ and $\widetilde{\beta}_n \cap \partial\widetilde{B}$ respectively converge to~$\widetilde{l} \cap \partial \widetilde{B}$ and~$\widetilde{l'} \cap \partial\widetilde{B}$ under the Hausdorff topology on $\partial \widetilde{B}$; thus, $\widetilde{\alpha}_n \cap \partial \widetilde{B}$ also cuts $\partial \widetilde{B}$ into two components such that the points $\widetilde{\beta}_n \cap \partial\widetilde{B}$ lie in different components for $n$ sufficiently large.
 This implies that some $\widetilde{\alpha}_n$ and~$\widetilde{\beta}_n$, hence $\alpha_n$ and $\beta_n$, intersect transversely, yielding a contradiction.

 Next, consider the foliation components $Y, Y' \subseteq S$ of $M_q(\theta)$ and $M_q(\theta')$ containing $l$ and~$l'$ respectively.
 Our strategy is to show that $Y$ and $Y'$ have disjoint interiors. This will imply that~$Y$ and $Y'$ are both proper subsurfaces,
 and so taking $\alpha$ and $\beta$ to be any saddle connection on their respective boundaries will satisfy the given requirements.
 Suppose for a contradiction that~$Y$ and $Y'$ have overlapping interiors. We argue that $l$ and $l'$ must intersect transversely, which is forbidden.
 If $Y$ and $Y'$ are both cylinders, then their core curves have non-zero geometric intersection number.
 If $Y$ and $Y'$ are not both cylinders, then without loss of generality, we can assume that~$l$ is dense in $Y$; in this case, $l$ must intersect~$Y'$ transversely.
 In both cases, we deduce that~$l$ and $l'$ intersect transversely as~claimed.
 \endproof

 Let $N_1(A)$ denote the 1--neighbourhood of a subset $A$ in $\G(S,q)$.

 \begin{lemma}[$1$--neighbourhoods of closed sets are closed] \label{lem:neighbourhood_of_closed_is_closed}
  Suppose that $A \subseteq \G(S,q)$ is closed in~$\RP$. Then $N_1(A)$ is also closed.
 \end{lemma}

 \proof
 Let $\theta_n \in N_1(A)$ be a sequence of slopes that converges to a slope $\theta$.
 Our goal is to show that~\hbox{$\theta \in N_1(A)$}.
 For each $\theta_n$, there exists some $\theta'_n \in A$ that is either adjacent to $\theta_n$ in~$\G(S,q)$ or equal to $\theta_n$.
 If $\theta_n = \theta'_n$ infinitely often then, since $A$ is closed, we have~$\theta \in A$ as required.
 If not, then by discarding finitely many terms, we may assume that $\theta_n \neq \theta'_n$ for all $n$.
 Therefore, for every~$n$, there exists a pair of disjoint saddle connections $\alpha_n$ and~$\beta_n$ with respective slopes~$\theta_n$ and~$\theta'_n$. By passing to a subsequence, we may assume that $\theta'_n$ converges to a slope $\theta'$.
 If $\theta = \theta'$ then we are done, so let us suppose otherwise.
 By the previous lemma, there exists a pair of disjoint saddle connections $\alpha$, $\beta$ with slopes $\theta$, $\theta'$.
 Since $A$ is closed, we have $\theta' \in A$, hence~$\theta \in N_1(A)$.
 \endproof

 We can now prove \cref{prop:balls_in_slope_graph_are_closed}.

 \begin{proof}[Proof of \cref{prop:balls_in_slope_graph_are_closed}]
  Let $\theta \in \G(S,q)$ and $r \geq 0$.
  As singletons are closed in $\RP$, an inductive application of \cref{lem:neighbourhood_of_closed_is_closed} shows that the ball $B(\theta, r)$
  is closed in $\RP$.
 \end{proof}

In particular, as $B(\theta, r)$ is closed and countable, its complement in $\RP$ is a disjoint union of open intervals.
We now show that every open interval arising as a component of ${\RP \setminus B(\theta, r)}$ for some $\theta\in \G(S,q)$ and $r \geq 0$ contains infinitely many slopes in $B(\theta, r+1)$.
For this, we require the following lemma.

\begin{lemma}[Converging from either side]
 Let $\alpha$ be a saddle connection with slope $\theta$.
 Then there exists a sequence of saddle connections $\alpha_n$, each disjoint from $\alpha$,
 with corresponding slopes~$\theta_n \neq \theta$ converging to $\theta$.
 Moreover, the $\alpha_n$ can be chosen so that $\theta_n$ converges to $\theta$ from either~side in~$\RP$.
\end{lemma}

 \proof
By applying a rotation, we may assume that $\alpha$ is vertical.
Since the set of slopes of saddle connections is dense in $\RP$ (see \cite[Theorem~2]{masur_86}), there exists a sequence $\beta_n \in \A(S,q)$ with positive slopes $\theta'_n$ converging to the vertical direction as $n \rightarrow \infty$.
By \cref{lem:straight_ladders}, for every~$n\in \mathbb{N}$, there exists a ladder path in $\A(S,q)$ from $\beta_n$ to $\alpha$ such that all saddle connections in this path except $\alpha$ have (positive) slope at least as steep as $\theta'_n$. In particular, there is a saddle connection $\alpha_n$ disjoint from~$\alpha$ with (positive) slope $\theta_n$ at least as steep as $\theta_n'$. Since $\theta_n' \rightarrow \theta$, we also have $\theta_n \rightarrow \theta$.
The proof for the negative slope case is analogous.
 \endproof

Now we can prove the promised proposition.

\begin{proof}[Proof of \cref{prop:infinitely_many_r+1}]
 By \cref{prop:balls_in_slope_graph_are_closed}, $B(\theta, r)$ is closed, hence any connected component~$I$ of $\RP \setminus B(\theta, r)$ is an open interval. Moreover, the endpoints of $I$ belong to $B(\theta, r)$; let $\theta_0$ be one endpoint. Applying the previous lemma, there exists a sequence $\theta_0 \neq \theta_n \in N_1(\theta_0)$ converging to~$\theta_0$ such that $\theta_n \in I$ for all $n$. In particular, $\theta_n \in B(\theta, r+1) \setminus B(\theta, r)$ for all $n$.
\end{proof}

Note that \cref{prop:infinitely_many_r+1} holds in more generality: For any closed $A \subseteq \G(S,q)$, a maximal open interval in $\RP \setminus A$ contains infinitely many slopes in $N_1(A)$.
The proof also shows directly that such an interval has no endpoint in common with any interval in~$\RP \setminus N_1(A)$.

\subsection{A partition of the graph of slopes}

In this subsection, we construct a partition of the graph of slopes using a family of intervals in~$\RP$.
For any (non-degenerate) interval $I \subseteq \RP$, let $G(I)$ be the set of vertices in $\G(S,q)$ whose slopes are contained in~$I$.
As the slopes of all saddle connections on $(S,q)$ are dense in~$\RP$, also $G(I)$ is dense as a subset of $I$.

 \begin{lemma}[Connected and infinite diameter] \label{lem:G_I_infinite_diameter}
  Let $I \subseteq \RP$ be an interval. Then the induced subgraph of $G(I)$ in $\G(S,q)$ is connected and has infinite diameter.
 \end{lemma}

 \proof
 Let $\alpha$ and $\beta$ be non-parallel saddle conections whose slopes both belong to $I$. By applying an $\SL$--deformation, we may assume that $\alpha$ is vertical and $\beta$ is horizontal. Then the (modified) interval $I$ contains all the positive slopes or all the negative slopes. Hence, by \cref{lem:straight_ladders}, all slopes appearing on one of the two ladder paths from $\alpha$ to $\beta$ lie in $I$. Therefore,~$G(I)$ spans a connected subgraph of $\G(S,q)$.

 Next, suppose for a contradiction that $G(I)$ has finite diameter.
 Then $G(I) \subseteq B(\theta, r)$ for some~$\theta \in \G(S,q)$ and $r \geq 0$.
 By \cref{prop:balls_in_slope_graph_are_closed}, $B(\theta, r)$ is a closed subset of $\RP$, and as it is also countable, $I\setminus G(I)\supseteq I\setminus B(\theta, r)$ contains a (non-degenerate) open interval. But this contradicts the fact that $G(I)$ is dense in $I$.
 \endproof

 We point out that the lemma above gives another proof that $\A(S,q)$ has infinite diameter.

 Next, we define a family of intervals using balls in the graph of slopes.
 Fix a basepoint ${\theta_0 \in \G(S,q)}$, and write $B(k) =  B(\theta_0, k)$.
 For each \hbox{$k \geq 0$}, let~$\I(k)$ be the set of connected components of $\RP \setminus B(k)$;
 by \cref{prop:balls_in_slope_graph_are_closed}, each such component is an open interval.
 Note that $\I(0)$ contains only the interval $\RP \setminus \set{\theta_0}$.
 The collection of intervals $\I = \sqcup_{k \geq 0}~ \I(k)$ forms a poset with respect to inclusion.

 Let $T$ be the Hasse diagram of $\I$, regarded as an (undirected) graph.
 That is, we have an edge between $I$ and~$I'$ if and only if $I \subsetneq I'$ and there exists no $I''$ such that $I \subsetneq I'' \subsetneq I'$ (or vice versa).

 \begin{lemma}[Hasse diagram is a tree] \label{lem:Hasse_tree}
  The graph $T$ is the regular (countably) infinite-valent tree $T_\infty$.
 \end{lemma}

 \proof
 For each $I \in \I$, the set $\set{I' \in \I \st I \subsetneq I'}$ is well-ordered.
 Furthermore, $\RP \setminus \set{\theta_0}$ is the unique maximal element of $\I$. Thus $T$ is a (rooted) tree.
 By \cref{prop:infinitely_many_r+1}, each $I \in \I$ has (countably) infinitely many edges.
 \endproof

 We are now ready to define a partition of the vertex set of $\G(S,q)$, where the equivalence classes are indexed by $\I$.
 For each $I\in \I$, define the associated \emph{slice} $Z(I)$ as follows.
 Set $Z(\RP \setminus \set{\theta_0}) := B(3)$.
 For each $k \geq 1$ and $I \in \I(k)$, set
 \[{Z(I) := G(I) \cap ( B(k+3) \setminus B(k+2) )};\]
 thus $\theta \in Z(I)$ if and only if $\theta \in I$ and $d(\theta_0, \theta) = k+3$.
 Our goal is to show that the slices indeed partition $G(S,q)$, are uniformly bounded, and give rise to a quotient graph isomorphic to $T$.

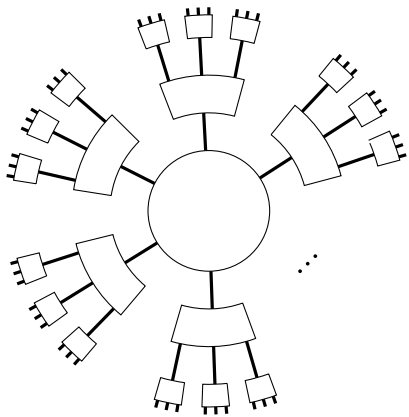
\begin{figure}
  \centering
  \begin{tikzpicture}
   \draw (0,0) circle (0.8cm);

   \draw (-70:1.3) arc (-70:-106:1.3) -- (-106:1.8) arc (-106:-70:1.8) -- (-70:1.3);
   \draw (-130:1.3) arc (-130:-166:1.3) -- (-166:1.8) arc (-166:-130:1.8) -- (-130:1.3);
   \draw (171:1.3) arc (171:135:1.3) -- (135:1.8) arc (135:171:1.8) -- (171:1.3);
   \draw (111:1.3) arc (111:75:1.3) -- (75:1.8) arc (75:111:1.8) -- (111:1.3);
   \draw (51:1.3) arc (51:15:1.3) -- (15:1.8) arc (15:51:1.8) -- (51:1.3);

   \foreach \x in {-88,-148,153,93,33}
   \draw[very thick] (\x:0.8) -- (\x:1.3);

   \draw (-70:2.3) arc (-70:-78:2.3) -- (-78:2.6) arc (-78:-70:2.6) -- (-70:2.3);
   \draw (-84:2.3) arc (-84:-92:2.3) -- (-92:2.6) arc (-92:-84:2.6) -- (-84:2.3);
   \draw (-98:2.3) arc (-98:-106:2.3) -- (-106:2.6) arc (-106:-98:2.6) -- (-98:2.3);

   \draw (-130:2.3) arc (-130:-138:2.3) -- (-138:2.6) arc (-138:-130:2.6) -- (-130:2.3);
   \draw (-144:2.3) arc (-144:-152:2.3) -- (-152:2.6) arc (-152:-144:2.6) -- (-144:2.3);
   \draw (-158:2.3) arc (-158:-166:2.3) -- (-166:2.6) arc (-166:-158:2.6) -- (-158:2.3);

   \draw (171:2.3) arc (171:163:2.3) -- (163:2.6) arc (163:171:2.6) -- (171:2.3);
   \draw (157:2.3) arc (157:149:2.3) -- (149:2.6) arc (149:157:2.6) -- (157:2.3);
   \draw (143:2.3) arc (143:135:2.3) -- (135:2.6) arc (135:143:2.6) -- (143:2.3);

   \draw (111:2.3) arc (111:103:2.3) -- (103:2.6) arc (103:111:2.6) -- (111:2.3);
   \draw (97:2.3) arc (97:89:2.3) -- (89:2.6) arc (89:97:2.6) -- (97:2.3);
   \draw (83:2.3) arc (83:75:2.3) -- (75:2.6) arc (75:83:2.6) -- (83:2.3);

   \draw (51:2.3) arc (51:43:2.3) -- (43:2.6) arc (43:51:2.6) -- (51:2.3);
   \draw (37:2.3) arc (37:29:2.3) -- (29:2.6) arc (29:37:2.6) -- (37:2.3);
   \draw (23:2.3) arc (23:15:2.3) -- (15:2.6) arc (15:24:2.6) -- (24:2.3);

   \foreach \x in {-74,-88,-102, -134,-148,-162, 167,153,139, 107,93,79, 47,33,19}
   \draw[very thick] (\x:1.8) -- (\x:2.3);

   \foreach \x in {-71,-74,-77, -85,-88,-91, -99,-102,-105,
   -131,-134,-137, -145,-148,-151, -159,-162,-165,
   170,167,164, 156,153,150, 142,139,136,
   110,107,104, 96,93,90, 82,79,76,
   50,47,44, 36,33,30, 22,19,16}
   \draw[very thick] (\x:2.6) -- (\x:2.7);

   \fill (1.2,-0.8) circle (0.7pt);
   \fill (1.3,-0.7) circle (0.7pt);
   \fill (1.4,-0.6) circle (0.7pt);

  \end{tikzpicture}
  \caption{Schematic of the partition of $\G(S,q)$ by slices $Z(I)$ for $I \in \I$.}
  \label{fig:partition_graph_of_slopes}
 \end{figure}

 \begin{lemma}[Slices partition $\G(S,q)$]
  The slices $\set{Z(I) \st I \in \I}$ form a partition of the vertex set of $\G(S,q)$.
 \end{lemma}

 \proof
 Let $\theta$ be a vertex of $\G(S,q)$.
 If $d(\theta, \theta_0) \leq 3$ then $Z(\RP \setminus \set{\theta_0}) = B(3)$ is the unique slice containing $\theta$.
 Now assume that $k\coloneqq d(\theta, \theta_0) - 3 \geq 1$. We have $\theta \in B(k+3) \setminus B(k+2)$ and in particular $\theta \notin B(k)$ which implies that there exists a unique interval $I \in \I(k)$ such that $\theta \in G(I)$.
 In either case, $\theta$ belongs to exactly one slice.
 \endproof

 Next, we identify the quotient graph induced by this partition. To do so, we need to characterise when two slices share a pair of adjacent vertices.

 \begin{lemma}[Slices and Hasse diagram]
  Let $I \neq I'$. Then $I$ and $I'$ are adjacent in $T$ if and only if there exists an edge in~$\mathcal{G}(S,q)$ connecting a pair of vertices $v\in Z(I)$ and $v'\in Z(I')$.
 \end{lemma}

 \proof
 Suppose that $I' \subsetneq I$ are adjacent in $T$. Then $I \in \I(k)$ and $I' \in \I(k+1)$ for a unique~$k \geq 0$.
 If $k = 0$ then every $v \in Z(I') = G(I') \cap ( B(4) \setminus B(3) )$ is at distance 4 from $\theta_0$, and so is adjacent to some vertex $v' \in B(3) = Z(I)$.
 Now assume that $k \geq 1$. Let $\theta \in \G(S,q)$ be an endpoint of $I'$ that is not an endpoint of $I$.
 Then $\theta \in G(I)$ and in particular $\theta \in B(k+1)\setminus B(k)$.
 Since $G(I')$ has infinite diameter in $\G(S,q)$ by \cref{lem:G_I_infinite_diameter}, we may choose some $\theta' \in G(I')$ outside $B(k+3)$. Then there exists a ladder path in the slope graph from $\theta$ to~$\theta'$ such that all intermediate vertices lie in $G(I') \subset G(I)$. In particular, this path has a pair of adjacent vertices~$v$ and~$v'$ at distance exactly~$k+3$ and $k+4$ from $\theta_0$ as desired.

 For the converse, suppose that $v\in Z(I)$ and $v'\in Z(I')$ are adjacent in $\G(S,q)$, where $I \in \I(k)$ and $I' \in \I(k')$.
 As $d(v,v') \leq 1$, we may assume, without loss of generality, that $k \leq k' \leq k+1$.
 Suppose, for a contradiction, that $I$ and $I'$ are not adjacent in $T$, and are thus disjoint as intervals.
 Then there exists a pair of slopes $\theta, \theta' \in B(k)$ such that $v$ and $v'$ are contained in different components of $\RP \setminus \{\theta, \theta'\}$. Since $B(k)$ is connected, we may choose $\theta$ and $\theta'$ to be adjacent.
 By \cref{lem:separating_slopes}, at least one of $v$ or $v'$ is within distance two of $\{\theta, \theta'\}$ in $\G(S,q)$,
 and is hence contained in $B(k+2)$. This contradicts the definition of the slices.
 \endproof

 Consequently, the map $\mathbf{I} \from \G(S,q) \to T$ sending all vertices in $Z(I)$ to $I$ realises the quotient map associated to the partition.
 Thus, the slices partition $\G(S,q)$ in a tree-like fashion as indicated in~\cref{fig:partition_graph_of_slopes}.
 Finally, in order to use \cref{lem:quotient_maps_quasi_isometries}, we must bound the diameter of each slice.

 \begin{lemma}[Slices have bounded diameter]
  For every $I \in \I$, the diameter of the slice $Z(I)$ in~$\G(S,q)$ is at most $17$.
 \end{lemma}

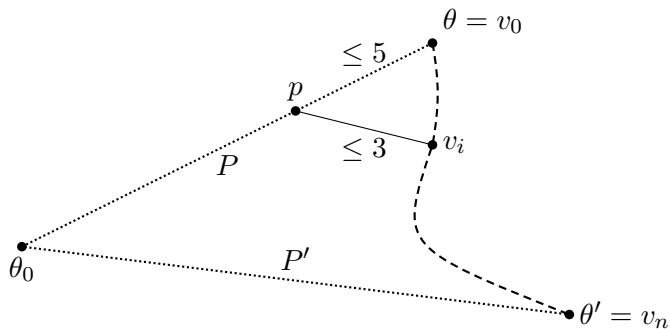
\begin{figure}
  \centering
  \begin{tikzpicture}[scale=0.9]
   \fill (0,0) circle (2pt) node[below]{$\theta_0$};

   \fill (6,3) circle (2pt) node[above right]{$\theta=v_0$};
   \fill (8,-1) circle (2pt) node[right]{$\theta'=v_n$};

   \draw[thick, densely dashed] (6,3) to[bend left=10] (6,1.5) .. controls (5.5,0) .. (8,-1);
   \fill (6,1.5) circle (2pt) node[right]{$v_i$};

   \draw[thick, densely dotted] (6,3) -- node[above]{$\leq 5$} (4,2) -- 
   (0,0) -- node[above]{$P'$} (8,-1);
   \fill (4,2) circle (2pt) node[above]{$p$};
   \draw (3,1.5) node[below]{$P$};

   \draw (6,1.5) -- node[below]{$\leq 3$} (4,2);
  \end{tikzpicture}
  \caption{Connecting $\theta$ and $\theta'$ by a ladder path (dashed), and by a concatenation of two geodesics (dotted).}
  \label{fig:slice_diameter}
 \end{figure}

 \proof
 Assume that $I \in \I(k)$.
 If $k \leq 4$ then $Z(I) \subseteq B(7)$, and so $\diam Z(I) \leq 14$.
 We may thus assume that $k \geq 5$. Let $\theta, \theta' \in Z(I) \subseteq G(I)$.
 Let $P$ and $P'$ be geodesics in~$\G(S,q)$ from~$\theta_0$ to~$\theta$ and $\theta'$ respectively (see \cref{fig:slice_diameter} for a sketch).
 By \cref{lem:straight_ladders} and projecting to $\G(S,q)$, there exists a ladder path $\theta = v_0, \ldots, v_n = \theta'$ in the slope graph with $v_i \in G(I)$ for all $i$.
 Note that $G(I)$ is disjoint from $B(k)$, and so $d(\theta_0, v_i) \geq k+1$ for each~$i$.
 Since the concatenation of $P$ with~$P'$ is also a path from $\theta$ to $\theta'$, by \cref{prop:slopes_ladder_paths_close_to_every_path} each~$v_i$ lies within distance~$3$ of some vertex~$p \in P \cup P'$. Without loss of generality, suppose~$p \in P$.
 Then
 \[d(\theta_0, p) \geq d(\theta_0, v_i) - d(v_i, p) \geq (k + 1) - 3 = k - 2.\]
 Since $P$ is a geodesic, we have
 \[d(\theta, p) = d(\theta, \theta_0) - d(\theta_0, p) \leq (k + 3) - (k - 2) = 5, \]
 hence $d(\theta, v_i) \leq d(\theta, p) + d(p, v_i) \leq 5 + 3 = 8$.
 Thus, every vertex along the ladder path lies in the $8$--neighbourhood of $\{\theta,\theta'\}$ in $\G(S,q)$,
 and so $d(\theta, \theta') \leq 8 + 1 + 8 = 17$.
 \endproof

It follows that the preimage of each slice in $\A(S,q)$ has diameter at most $35$. However, we shall briefly repeat the above proof for the saddle connection graph to obtain stronger bounds.

 \begin{lemma}[Preimages of slices have bounded diameter]
  For every $I \in \I$, the diameter of the preimage of the slice $Z(I)$ under $\Theta$ is at most $29$ in $\A(S,q)$.
 \end{lemma}

 \proof
 Assume that $I \in \I(k)$.
 If $k \leq 4$ then $\diam \Theta^{-1}(Z(I)) \leq 2 \cdot \diam Z(I) + 1 = 2 \cdot 14 + 1 = 29$.
 Now assume $k \geq 5$.
 Let $\alpha, \alpha' \in \A(S,q)$ be saddle connections with slopes $\theta, \theta' \in Z(I)$.
 Let~$\hat P, \hat P'$ be paths in $\A(S,q)$ that descend to the geodesics $P,P'$ in $\G(S,q)$ as in the previous lemma.
 We may similarly choose a ladder path $\alpha = a_0, \ldots, a_n = \alpha'$ in $\A(S,q)$ that descends to the ladder path between $\theta$ and $\theta'$ as given above. By \cref{prop:ladder_paths_close_to_every_path}, each $a_i$ lies within distance~$3$ of some $\hat p \in \hat P \cup \hat P'$.
 We argue as before (and assume $\hat p \in \hat P$): As $d(\Theta(\hat p), \Theta(a_i)) \leq 3$, we get $d(\theta, \Theta(\hat p)) \leq 5$ in $\G(S,q)$. With this, we see that $d(\alpha, \hat p) \leq 2 \cdot d(\theta, \Theta(\hat p)) + 1 \leq 11$, hence $d(\alpha, a_i) \leq d(\alpha, \hat p) + d(\hat p, a_i) \leq 11 + 3 = 14$. Therefore, $d(\alpha, \alpha') \leq 14 + 1 + 14 = 29$.
 \endproof

 We can now use \cref{lem:quotient_maps_quasi_isometries} to obtain the following.

 \begin{corollary}[$\A(S,q)$, $\G(S,q)$, and $T$ are quasi-isometric]\label{cor:q_i}
  The maps $\mathbf{I} \from \G(S,q) \to T$ and $\mathbf{I}\circ\Theta \from \A(S,q) \to T$ are respectively a $(18,17)$--quasi-isometry and a $(30,29)$--quasi-isometry. \qed
 \end{corollary}

 This result combined with \cref{lem:Hasse_tree} completes the proof of \cref{thm:main}.

\section{Boundary structure}

In this section, we characterise the Gromov boundary of the saddle connection graph~$\A(S,q)$.
By \cref{prop:MT_2}, there is an isometric embedding $\jmath \from \A(S,q) \to \AC(S,\ZZ)$, which induces an embedding $\partial \jmath \from \partial\A(S,q)\to \partial\AC(S,\ZZ)  \cong \partial \C(S,\ZZ)$ (see \cref{prop:Pap}).
By \cref{thm:PO}, there is a natural identification of $\partial\C(S,\ZZ)$ with the space of arational topological foliations~$\AF(S,\ZZ)$.
Our goal is to characterise the foliations that appear in the image of $\partial \jmath$.
We shall show that it is exactly those that arise from straight foliations on $(S,q)$ that contain no saddle connections as~leaves.

Throughout this section, we shall consider three sequences simultaneously:
\begin{itemize}
 \item a sequence $\alpha_n \in \A(S,q)$ of saddle connections,
 \item the corresponding slopes $\theta_n = \Theta(\alpha_n) \in \G(S,q)$, and
 \item the associated intervals $I_n = \mathbf{I}(\theta_n) \in T$.
\end{itemize}
As the saddle connection graph $\A(S,q)$, the graph of slopes $\G(S,q)$, and the Hasse diagram $T$ (for a collection of intervals $\I$ with a fixed basepoint $\theta_0$) are all quasi-isometric by \cref{cor:q_i}, each of the three sequences is admissible (in the appropriate graph) if and only if the other two are.
Let us now characterise admissibility in each of these cases.

\pagebreak

\begin{lemma}[Admissibility for intervals]
 Let $I_n$ be a sequence of intervals in $T$. The following are equivalent.
 \begin{enumerate}
  \item The sequence $I_n$ is admissible in $T$,
  \item there exists a geodesic ray $\RP \setminus \set{\theta_0} = J_0, J_1, J_2, \ldots$ in $T$ such that for each~$k \geq 0$,
  the unique geodesic from $\RP \setminus \set{\theta_0}$ to $I_n$ in $T$ passes through $J_k$ for all $n$ sufficiently large, and
  \item there exists an arational slope $\theta\in\RP$ such that for any open neighbourhood $U \subseteq \RP$ of~$\theta$, we have $I_n \subseteq U$ for all $n$ sufficiently large.
 \end{enumerate}
 Furthermore, when the above hold, $\theta$ is unique and we have $\cap_k J_k = \set{\theta}$ in $\RP$.
\end{lemma}

\proof
``(i) $\Leftrightarrow$ (ii)'': The sequence $I_n$ is admissible if and only if for all $k \geq 0$, there exists some~$N_k\geq 0$ such that the Gromov product $(I_n | I_m)_{\RP \setminus \set{\theta_0}} \geq k$ for all $m,n \ge N_k$.
Since~$T$ is a tree, this is equivalent to saying that for all $k \geq 0$, there exists some $J_k\in T$ at distance~$k$ from~$\RP \setminus \set{\theta_0}$ such that the unique geodesic from $\RP \setminus \set{\theta_0}$ to $I_n$ in $T$ passes through $J_k$ for all~$n$ sufficiently large. Note that when this holds, the sequence $J_k$ forms a geodesic ray.
Thus, the first two statements are equivalent.

``(ii) $\Rightarrow$ (iii)'':
By construction of $T$, the sequence $J_k$ forms a descending nested sequence of open intervals.
Therefore, the intersection of their closures $\cap_k \overline{J_k}$ is non-empty. By the remark after \cref{prop:infinitely_many_r+1}, $J_{k+1}$ has no endpoints in common with $J_k$ for all $k$. Therefore, any slope in~$\cap_k \overline{J_k}$ cannot be an endpoint of $\overline{J_k}$ for all $k$. Thus, $\cap_k J_k = \cap_k \overline{J_k} \neq \emptyset$. Now, $\cap_k J_k$ is a connected subset of~$\RP \setminus \G(S,q)$. Since $\G(S,q)$ is dense in $\RP$, it follows that $\cap_k J_k = \set{\theta}$ for some $\theta \in \RP \setminus \G(S,q)$, i.e.\ some arational slope $\theta$. Let~$U\subseteq \RP$ be any open neighbourhood of~$\theta$. Since the sequence of balls $B(\theta_0, k)$ exhausts the dense set~$\G(S,q)$, there exists some $k$ such that~$\theta$ lies in an interval from $\I(k)$ having both endpoints in $U$; in fact, this interval must be $J_k$. Since the unique geodesic from $\RP \setminus \set{\theta_0}$ to $I_n$ in~$T$ passes through~$J_k$, we see that $I_n \subseteq J_k \subseteq U$ for all~$n$ sufficiently large.

``(iii) $\Rightarrow$ (ii)'': Assume that $\theta$ is an arational slope satisfying Property (iii). For each~$k \geq 0$, let~$J_k$ be the unique (open) interval from $\I(k)$ containing $\theta$. Then, by assumption, $I_n \subseteq J_k$ for all $n$ sufficiently large. In particular, this implies that the geodesic from $\RP \setminus \set{\theta_0}$ to~$I_n$ in $T$ passes through $J_k$ for all $n$ sufficiently large.
\endproof

Note that each point on the Gromov boundary of $T$ is represented by a unique geodesic ray starting from $\RP \setminus \set{\theta_0}$. Thus, two admissible sequences $I_n$ and $I'_n$ in~$T$ are equivalent if and only if they satisfy Property (iii) above for a common arational slope~$\theta$.

\begin{lemma}[Admissibility for slopes]
  A sequence of slopes $\theta_n$ is admissible in $\G(S,q)$ if and only if $\theta_n$ converges to an arational slope $\theta \in \RP$. Moreover, when this holds then~$\theta$ is the unique slope contained in infinitely many of the intervals $I_n = \mathbf{I}(\theta_n)$.
\end{lemma}

\proof
 Assume that $\theta_n$ is an admissible sequence in $\G(S,q)$. Then the associated sequence of intervals~$I_n=\mathbf{I}(\theta_n)$ is also admissible in $T$. Therefore, we may define the (unique) sequence~$J_k$ as in the previous lemma satisfying $\cap_k J_k = \set{\theta}$ for some arational slope $\theta$.
 Since for each $k\geq 0$, we have $\theta_n \in I_n \subseteq J_k$ for all $n$ sufficiently large, we deduce that $\theta_n$ converges to $\theta$.

 For the converse, for each $k \geq 0$, let $J_k$ be the unique interval from $\I(k)$ containing $\theta$. Since~$\theta_n$ converges to~$\theta$, there exists $N_k \geq 0$ such that $\theta_n \in J_k$ for all $n \geq N_k$. By choosing $N_k \geq k$, we can also ensure that $\theta_n \in I_n \subseteq J_k$ for all $n \geq N_k$. Then by the previous lemma, we deduce that~$I_n$, hence $\theta_n$, are admissible.
\endproof

Similarly to admissible sequences in $T$, two admissible sequences of slopes are equivalent if and only if their limits in $\RP$ coincide.

\begin{proposition}[Admissibility for saddle connections]
 Let $\alpha_n$ be a sequence of saddle connections on $(S,q)$.
 The following are equivalent:
 \begin{enumerate}
  \item The sequence $\alpha_n$ is admissible in $\A(S,q)$ and
  \item the corresponding slopes $\theta_n = \Theta(\alpha_n)$ converge to an arational slope $\theta$.
 \end{enumerate}
 
 Moreover, when the above hold, the sequence $\alpha_n$ converges to an arational topological foliation $F_\infty \in \partial\A(S,q)\subset\partial\AC(S,\ZZ)$ and the straight foliation on $(S,q)$ with slope $\theta$ is topologically equivalent to $F_\infty$.
\end{proposition}

\proof
The equivalence of (i) and (ii) follows from the previous lemma and \cref{cor:q_i}.
That~$\alpha_n$ converges to an arational topological foliation follows directly from \cref{thm:PO}.

Assume now that the sequence $\alpha_n$ is admissible and hence $\theta_n$ converges to an arational slope~$\theta$.
Let $(F_n, \mu_n) = M_q(\theta_n)$ and $(F, \mu) = M_q(\theta)$.
Then $F$ contains no saddle connections and by \cref{lem:arational}, $F$ is an arational foliation.
Our goal is to prove that $F$ and $F_\infty$ are topologically equivalent.
Let~$\gamma_n \in \C(S,\ZZ)$ be a boundary component of a small regular neighbourhood of $\alpha_n$ on~$(S,\ZZ)$.
Observe that
\[i(\gamma_n, \alpha_n) = i(\gamma_n, (F_n, \mu_n)) = 0\]
for all $n$.
Thus, the sequences $\alpha_n$ and $\gamma_n$ fellow-travel in $\AC(S,\ZZ)$,
and so they both converge to the same arational foliation $F_\infty \in \partial \AC(S,\ZZ)$.
By passing to a sub\-sequence and appealing to compactness of $\PMF(S,\ZZ)$,
there exists a measured foliation $(F', \mu') \in \MF(S,\ZZ)$ with \mbox{$\gamma_n \rightarrow (F', \mu')\in\PMF(S,\ZZ)$}. In particular, there exist constants $c_n > 0$ such that 
\[c_n\gamma_n \rightarrow (F', \mu')\in\MF(S,\ZZ)\]
as $n \rightarrow \infty$.
Therefore, we have
\[i((F', \mu'),(F,\mu))=\lim_{n\to\infty}
c_n\cdot i(\gamma_n,(F_n,\mu_n))=0.\]
On the one hand, by \cref{lem:arational_2} and since $F$ is arational, $F'$ and $F$ are topologically equivalent. On the other hand, by \cref{thm:PO}, $F'$ is also topologically equivalent to~$F_\infty$. Hence, $F$ is topologically equivalent to $F_\infty$ as required.
\endproof

 Thus, the slope map $\Theta \from \A(S,q) \to \RP$ has an extension
 \[\partial \Theta \from \partial \A(S,q) \to \RP \setminus \mathrm{im}(\Theta) \subset \RP\]
 defined by assigning each (equivalence class of an) admissible sequence its limiting slope.
 Since the slopes of saddle connections are dense in $\RP$, the map
 $\partial \Theta$ is a bijection (but not a homeomorphism).

\bibliographystyle{amsalpha}
\bibliography{literature_boundary_saddle_connection_graph}

\end{document}